\documentclass[preprint,1p,sort&compress]{elsarticle}
\usepackage{lineno,hyperref}

\usepackage{graphicx} 
\usepackage{amsmath} 
\usepackage{amssymb}  
\usepackage{epstopdf}
\usepackage{multirow}
\usepackage{amsthm}
\usepackage{subfigure}
\usepackage{colortbl}
\usepackage{array}
\usepackage{eufrak}
\usepackage{setspace}
\usepackage{soul}
\usepackage{booktabs}
\usepackage{enumerate}
\usepackage{multicol}
\usepackage{makecell}

\journal{Journal of Process Control}
\newdefinition{rmk}{Remark}
\newproof{pf}{Proof}
\newproof{pot}{Proof of Theorem \ref{thm2}}

\begin{document}
\begin{spacing}{1.45}
\begin{frontmatter}
\title{A Data-Driven Robust Optimization Approach to Scenario-Based Stochastic Model Predictive Control}
\author[Tsinghua]{Chao Shang}\ead{c-shang@tsinghua.edu.cn}
\author[Cornell]{Fengqi You\corref{cor1}}\ead{fengqi.you@cornell.edu}
\address[Tsinghua]{Beijing National Research Center for Information Science and Technology, and Department of Automation, Tsinghua University, Beijing 100084, China}
\address[Cornell]{Robert Frederick Smith School of Chemical and Biomolecular Engineering, Cornell University, Ithaca, New York 14853, USA}
\cortext[cor1]{To whom all correspondence should be addressed: Tel:
+1 (607)255-1162; Fax: +1 (607)255-9166; E-mail:
fengqi.you@cornell.edu}
\begin{abstract}
Stochastic model predictive control (SMPC) has been a promising solution to complex control problems under uncertain disturbances. However, traditional SMPC approaches either require exact knowledge of probabilistic distributions, or rely on massive scenarios that are generated to represent uncertainties. In this paper, a novel scenario-based SMPC approach is proposed by actively learning a data-driven uncertainty set from available data with machine learning techniques. A systematical procedure is then proposed to further calibrate the uncertainty set, which gives appropriate probabilistic guarantee. The resulting data-driven uncertainty set is more compact than traditional norm-based sets, and can help reducing conservatism of control actions. Meanwhile, the proposed method requires less data samples than traditional scenario-based SMPC approaches, thereby enhancing the practicability of SMPC. Finally the optimal control problem is cast as a single-stage robust optimization problem, which can be solved efficiently by deriving the robust counterpart problem. The feasibility and stability issue is also discussed in detail. The efficacy of the proposed approach is demonstrated through a two-mass-spring system and a building energy control problem under uncertain disturbances.
\end{abstract}
\begin{keyword}
 Stochastic model predictive control
 \sep Chance constraints
 \sep Scenario programs
 \sep Robust model predictive control
 \sep Machine learning
\end{keyword}
\end{frontmatter}
\section{Introduction}
\label{sec:introduction}
Model predictive control (MPC) has become a commonly accepted technology in various industrial scenarios thanks to its wide applicability and effectiveness in addressing complex optimal control problems subject to input and state constraints \cite{maciejowski2002predictive,qin2003survey,rawlings2009model,chu2015model}. In the context of MPC, process behaviors are characterized by various dynamic models, such as step response models, transfer function models, and state-space models, based on which the dependence of system evolutions in the future on control actions can be explicitly predicted \cite{morari1999model,liu2009distributed,pourkargar2017impact}. Based on the idea of setting the initial state to be the current state of the process, the optimal control sequence can be derived by solving an optimization problem at every control instant, and the first control action is implemented only. This procedure is repeated at the next control instant, which results in a receding horizon control (RHC) scheme.

In practical scenarios, the system to be controlled is usually prone to different kinds of uncertain disturbances, which tend to make process states deviate from their nominal trajectory \cite{christofides2005control}. Unfortunately, this issue is not accounted for by deterministic MPC algorithms, and hence degradations of control performance or violations of system constraints may frequently occur \cite{kadam2007dynamic}. To immunize control actions against unknown disturbances, robust MPC (RMPC) and stochastic MPC (SMPC) algorithms have gained increasing attentions \cite{pannocchia2003robust,petersen2014robust,gonccalves2016robust}. In the context of RMPC, the uncertainty is described by a bounded set, based on which the worst-case control objective or constraint satisfaction is pursued. However, over-conservatism is inevitably induced as a downside of RMPC \cite{saltik2018outlook}. To alleviate the conservatism of RMPC, chance constraints are adopted in SMPC, which resort to probabilistic descriptions of disturbances \cite{farina2016stochastic}. However, the exact knowledge of distributional information is commonly unavailable in practice, and complicated computations of multi-dimensional integrals are involved, thereby rendering optimization problems difficult to solve.

To tackle the realistic implementation difficulty of SMPC, the scenario-based approximation strategy has been established in recent years because of its freedom from specific assumptions on distributions of disturbances \cite{calafiore2006scenario,campi2008exact,calafiore2010random}. In the scenario-based approach, chance constraints are replaced with a series of hard constraints based on sample realizations of uncertainties, thereby giving rise to an approximated solution. The major attractiveness is that the strict assumption on the probability distribution can be alleviated, and scenario-based programs are easier to handle. Applications of scenario-based SMPC can be found in smart grid operation \cite{zhu2014decomposed} and automotive control \cite{di2014stochastic,moser2018flexible}. Despite of its popularity, the scenario-based approach is subject to the possibility that the approximated solution fails to satisfy the chance constraints due to the randomness of data sampling. To ensure chance constraints to hold with predefined confidence levels, sufficient data samples must be collected \cite{calafiore2005uncertain,campi2008exact,calafiore2010random}, which lead to some obvious deficiencies. On one hand, the theoretical number of data samples that are required to safeguard the confidence level tends to be excessively large, resulting in massive hard constraints in the optimization problem. It not only poses considerable challenges for online computations, but also induces excessive conservatism of the solution. On the other hand, one may not be allowed to collect sufficient data samples due to some realistic limitations.

To address these practical challenges in classic SMPC techniques, we propose in this work a systematic data-driven robust optimization approach to SMPC. It seamlessly integrates machine learning into the optimization-based paradigm, which has received increasing attentions recently \cite{ning2017data,ning2018adaptive,shang2018distributionally}. Instead of solving the scenario-based program directly, we first learn a polytopic high-density region from available scenarios with support vector clustering (SVC), which has been proposed as a data-driven uncertainty set tailored to robust optimization \cite{shang2017data}. Then a novel calibration scheme is developed to further adjust the size of the uncertainty set based on an additional independent dataset, which establishes desirable probabilistic guarantee for SMPC. Finally, an alternative RMPC problem is solved based on the calibrated SVC-based uncertainty set, which gives an approximated solution to the original SMPC problem. In comparison with classic SMPC techniques, the proposed approach commonly requires a much lower sampling complexity to reach the probabilistic guarantee, thereby reducing the conservatism and enhancing the applicability of SMPC. Another advantage is that, since the data-driven uncertainty set can be expressed as a series of linear inequalities, the induced RMPC problem can be efficiently solved online with robust optimization techniques, thereby leading to a significant reduction of computational burden. Meanwhile, the feasibility and stability issue of the proposed control scheme is addressed rigorously. In the case of unbounded disturbance, system states may be driven towards infinity, resulting in infeasible optimization problems. On this occasion, an additional backup controller must be introduced to ensure feasibility for unbounded initial states. Then we establish system stability in a mean-square sense under mild assumptions on the performance of the backup controller and the properties of disturbances. The efficacy of the proposed SMPC method is demonstrated on a two-mass-spring system and a building energy control problem.


The layout of this article is organized as follows. In Section 2, we make clear the background on SMPC with chance constraints and the active uncertainty set construction approach based on SVC. The proposed methodology for SMPC is described in Section 3, while numerical results on a simulated two-mass-spring system are reported in Section 4. Section 5 presents a building control example under weather prediction uncertainty. Concluding remarks are given in Section 6.

\textbf{Notations and definitions}. Hereafter, $\mathbb{N}_0 = \{0, 1, \cdots \}$ denotes the set of non-negative natural numbers. $Q \succ (\succeq)\ 0$ denotes positive (semi-)definiteness of matrix $Q$. The $p$-norm of a matrix is denoted by $\|\cdot\|_p$, while $\|\cdot\|$ represents the Euclidean norm of a vector by convention. Given $Q \succeq 0$, $\|x\|_Q^2$ denotes the quadratic function $x^{\rm T}Qx$. The Kronecker product operator is denoted by $\otimes$. $I_m$ denotes the $m$-dimensional identity matrix, and $1_m$ denotes the $m$-dimensional vector with all elements being ones. The trace operator on a square matrix is denoted by ${\rm Tr}\{\cdot\}$.

\section{Preliminaries}
\subsection{Problem Formulation of SMPC with Chance Constraints}
In this work, we consider the predictive control problem of the following linear discrete-time time-invariant system subject to additive external uncertainties:
\begin{equation}
\begin{split}
x_{t+1} = Ax_t + B_u u_t + B_w w_t
\end{split}
\label{eq:statespace}
\end{equation}
where $x_t \in \mathbb{R}^{n_x}$ and $u_t \in \mathbb{R}^{n_u}$ are the state vector and input vector at stage $t$, respectively. $w_t \in \mathbb{R}^{n_w}$ denotes the uncertain disturbance vector. It is assumed that at each sampling instant we have access to the measurement of state $x_t$, and $(A,B_u)$ is stabilizable. The disturbance $\{w_t\}$ is assumed to have zero mean $\mathbb{E} \{ w_t \} = 0$. Different from standard assumptions in SMPC, here disturbances $\{ w_t \}$ are not required to be independently and identically distributed (i.i.d.). Instead, a data-driven approach is to be carried out in this paper to handle inherent temporal structures within the non-i.i.d. disturbance sequence. Given a control horizon $H$, the evolution of the system can be described in a matrix form:
\begin{equation}
\bold{x} = \bold{A}x_0 + \bold{B}_u\bold{u} + \bold{B}_w\bold{w}.
\label{eq:system}
\end{equation}
where the vector $\bold{x} = [x_1^{\rm T}\ x_2^{\rm T}\ \cdots \ x_H^{\rm T} ]^{\rm T} \in \mathbb{R}^{Hn_x}$ is obtained by stacking all state vectors within the control horizon $1 \le t \le H$, and $\bold{u} = [u_0^{\rm T}\ u_1^{\rm T}\ \cdots \ u_{H-1}^{\rm T} ]^{\rm T} \in \mathbb{R}^{Hn_u}$, $\bold{w} = [w_0^{\rm T}\ w_1^{\rm T}\ \cdots \ w_{H-1}^{\rm T} ]^{\rm T} \in \mathbb{R}^{Hn_w}$. System matrices $\bold{A}$, $\bold{B}_u$, and $\bold{B}_w$ are given by:
\begin{equation}
\begin{split}
\bold{A} = & \left [ \begin{array}{c}
A \\
A^2 \\
\vdots \\
A^H
\end{array}
\right ],\ \bold{B}_u = \left [ \begin{array}{c c c c}
B_u & 0 & \cdots & 0 \\
AB_u & B_u & \cdots & 0 \\
\vdots & \vdots & \ddots & \vdots \\
A^{H-1}B_u & A^{H-2}B_u & \cdots & B_u
\end{array}
\right ], \\
& \ \ \ \ \bold{B}_w = \left [ \begin{array}{c c c c}
B_w & 0 & \cdots & 0 \\
AB_w & B_w & \cdots & 0 \\
\vdots & \vdots & \ddots & \vdots \\
A^{H-1}B_w & A^{H-2}B_w & \cdots & B_w
\end{array}
\right ].
\end{split}
\end{equation}
Besides, we assume that in each stage $t$ the system states are subject to the following joint chance constraints:
\begin{equation}
\mathbb{P} \left \{ x_t \in F_x \right \} \ge 1 - \epsilon_t,\ t = 1,\cdots,H
\label{eq:cc}
\end{equation}
where $F_x = \{ x | Fx \le f \} \subseteq \mathbb{R}^{n_x}$ is a polytope and $F \in \mathbb{R}^{m \times n_x}$, $f \in \mathbb{R}^{m}$ are coefficients. $\epsilon_t \in [0,1)$ is the pre-specified maximal probability of constraint violation. One can also express (\ref{eq:cc}) in terms of $\bold{x}$:
\begin{equation}
\mathbb{P} \left \{ \bold{F}_t \bold{x} \le \bold{f}_t \right \} \ge 1 - \epsilon_t,\ t = 1,\cdots,H
\label{eq:cc1}
\end{equation}
where $\bold{F}_t$ and $\bold{f}_t$ denote the $t$th block row of $\bold{F} = I_H \otimes F$ and $\bold{f} = 1_H \otimes f$. Different from deterministic constraints, chance constraints allow partial violations that help alleviating the over-pessimism or infeasibility of hard constraints, especially when the uncertainty is unbounded. Meanwhile, hard input constraints are imposed:
\begin{equation}
u_t \in F_u = \left \{ u|G u \le g \right \} \subseteq \mathbb{R}^{n_u},\ t = 1,\cdots,H
\end{equation}
which can also compactly rewritten as:
\begin{equation}
\bold{Gu} \le \bold{g},
\label{eq:inputconstraints}
\end{equation}
where $\bold{G} = I_{n_u} \otimes G$ and $\bold{g} = 1_{n_u} \otimes g$. Different from $x_t$, manipulated variables $u_t$ are always subject to actuator limitations in practice, and hence using hard constraints bears realistic rationality. Moreover, admissible inputs are assumed to be bounded such that $\sup_{u \in F_u} \|u\| < \infty$.

A variety of cost functions can be chosen as the optimization objective. The most frequently used one in the SMPC literature is given by:
\begin{equation}
J = \mathbb{E} \left \{ \sum_{t=0}^{H-1} l(x_t,u_t) + l_f(x_H) \right \},
\label{eq:obj1}
\end{equation}
where $l(x_t,u_t)$ denotes the stage cost and $l_f(x_H)$ stands for the terminal cost. A popular choice is to adopt quadratic functions $l(x_t,u_t) \triangleq \|x_t\|^2_Q + \|u_t\|^2_R$ and $l_f(x_H) \triangleq \|x_H\|^2_{Q_f}$, where $Q \succeq 0,\ R \succ 0,\ Q_f \succeq 0$ with appropriate dimensions. Other cost functions defined based on nominal trajectories or sample-based averaging can also be applied \cite{mayne2016robust}. Given prediction horizon $H$, the optimization problem of SMPC to be considered in this paper can be formally cast as follows:
\begin{equation*}
\begin{split}
\min_\bold{u} &\ J {\rm\ defined\ in\ (\ref{eq:obj1})} \\
{\rm s.t.} &\ {\rm System\ dynamics\ (\ref{eq:system})} \\
           &\ {\rm State\ constraints\ (\ref{eq:cc1})} \\
           &\ {\rm Input\ constraints\ (\ref{eq:inputconstraints})} \\
\end{split}
\tag{SMPC}
\end{equation*}

\subsection{Disturbance Feedback Parameterization}
In principle, the future control action $u_t\ (t = 1,\cdots,H)$ shall be determined after the uncertain disturbances $\{w_0,\cdots,w_{t-1}\}\ (t = 1,\cdots,H)$ are observed, thereby depending on future realizations of disturbances. To taking such causal dependence into account, the disturbance feedback (DF) \textit{control policy} $u_t = \pi_t(w_1,\cdots,w_{t-1}) \ (t = 1,\cdots,H)$ is commonly used, where the causality here indicates that $u_t\ (t = 1,\cdots,H)$ is unrelated to future disturbances $w_t,w_{t+1},\cdots\ (t = 1,\cdots,H)$ that follow. A commonly adopted policy is the affine DF policy \cite{mayne2005robust,goulart2006optimization}:
\begin{equation}
u_t = h_t + \sum_{j=1}^{t-1} M_{t,j}w_j,\ t = 1,\cdots,H
\label{eq:ADF}
\end{equation}
where the matrix $M_{t,j}$ and the vector $h_t$ become decision variables hereafter. The affine DF policy is similar to the affine decision rule in multi-stage adaptive robust optimization \cite{ben2004adjustable}. By considering feedback policies for all stages as a whole, (\ref{eq:ADF}) can be further expressed as:
\begin{equation}
\bold{u} = \bold{h} + \bold{Mw}
\end{equation}
where
\begin{equation}
\bold{h} = \left [ \begin{array}{c}
h_0 \\
h_1 \\
\vdots \\
h_{H-1}
\end{array}
\right ],\ \bold{M} = \left [ \begin{array}{c c c c}
0 & 0 & \cdots & 0 \\
M_{1,0} & 0 & \cdots & 0 \\
\vdots & \vdots & \ddots & \vdots \\
M_{H-1,0} & M_{H-1,H-2} & \cdots & 0
\end{array}
\right ]
\end{equation}
are the constant coefficient vector and the lower triangular matrix, respectively, which serve as decision variables of the optimal control problem. Unfortunately, if the disturbance $\{ w_t \}$ has an unbounded support, system inputs will be unbounded, thereby inevitably resulting in violations of hard constraints (\ref{eq:inputconstraints}). To address this issue, Ref. \cite{hokayem2009stochastic} developed the so-called saturated DF policy:
\begin{equation}
\bold{u} = \bold{h} + \bold{M}\boldsymbol{\varphi}(\bold{w}),
\end{equation}
where for $\bold{w} = [w(0),\cdots,w(H-1)]^{\rm T}$, $\boldsymbol{\varphi}(\bold{w}) = [\varphi(w(0)),\cdots,\varphi(w(H-1))]^{\rm T}$ with $\varphi(\cdot)$ being a bounded saturation function that fulfills $|\varphi(\cdot)| \le 1$. Common options for the saturation function include the hyperbolic tangent function $\varphi(r) \triangleq \tanh(r)$ and the piecewise linear function $\varphi(r) \triangleq \min\{1,\max\{-1,r\}\}$. Under saturated DF, input constraints can be enforced in a robust manner \cite{paulson2017stochastic}:
\begin{equation}
\sup_{|\varphi(\cdot)| \le 1} \bold{G}\left [ \bold{h} + \bold{M}\boldsymbol{\varphi}(\bold{w}) \right ] \le \bold{g}.
\label{eq:robustinputconstraints}
\end{equation}
At the same time, the system dynamics model (\ref{eq:system}) becomes:
\begin{equation}
\bold{x} = \bold{A}x_0 + \bold{B}_u\bold{M}\boldsymbol{\varphi}(\bold{w}) + \bold{B}_w\bold{w} + \bold{B}_u\bold{h}.
\label{eq:allx}
\end{equation}
Substituting (\ref{eq:allx}) into (\ref{eq:cc1}), the chance constraint in the $t$th stage can be expressed as:
\begin{equation}
\mathbb{P} \left \{ \bold{F}_t\left (
\bold{A}x_0 + \bold{B}_u\bold{M}\boldsymbol{\varphi}(\bold{w}) + \bold{B}_w\bold{w} + \bold{B}_u\bold{h}
\right ) \le \bold{f}_t \right \} \ge 1 - \epsilon_t.
\label{eq:CC}
\end{equation}

With quadratic functions adopted as cost functions, and under the assumption that $\mathbb{E} \{ \boldsymbol{\varphi}(\bold{w}) \} = \bold{0}$,\footnote{This assumption can be always satisfied by subtracting the mean of $\boldsymbol{\varphi}(\bold{w})$.} the objective function (\ref{eq:obj1}) can be written as:
\begin{equation}
\begin{split}
& J(\bold{M},\bold{h}) \\
= &\ \mathbb{E} \left \{ \sum_{t=0}^{H-1}(x_t^{\rm T}Qx_t + u_t^{\rm T}Ru_t) + x_H^{\rm T}Q_fx_H \right \} \\
= &\ \mathbb{E} \left \{ \bold{x}^{\rm T} \bold{Q} \bold{x} + \bold{u}^{\rm T} \bold{R} \bold{u} \right \} \\
= &\ \bold{h}^{\rm T} \left ( \bold{R} + \bold{B}_u^{\rm T}\bold{QB}_u \right ) \bold{h} + {\rm Tr} \left \{ (\bold{R} + \bold{B}_u^{\rm T}\bold{QB}_u) \bold{M} \mathbb{E} \left \{ \boldsymbol{\varphi}(\bold{w})\boldsymbol{\varphi}(\bold{w})^{\rm T} \right \} \bold{M}^{\rm T} \right \} \\
&\ + 2\bold{h}^{\rm T}\bold{B}^{\rm T}_u \bold{QA}x + 2{\rm Tr} \left \{ \bold{M}^{\rm T}\bold{B}^{\rm T}_u\bold{QB}_w \mathbb{E} \left \{ \bold{w}\boldsymbol{\varphi}(\bold{w})^{\rm T} \right \} \right \} + C
\label{eq:J1}
\end{split}
\end{equation}
where $C = x^{\rm T}_0\bold{A}^{\rm T}\bold{QA}x_0 + {\rm Tr} \left \{ \bold{B}^{\rm T}_w\bold{QB}_w \mathbb{E} \left \{ \bold{w}\bold{w}^{\rm T} \right \} \right \}$ is a constant independent from $\bold{M}$ and $\bold{h}$.\footnote{The second-order moment information matrices $\mathbb{E} \left \{ \boldsymbol{\varphi}(\bold{w})\boldsymbol{\varphi}(\bold{w})^{\rm T} \right \}$ and $\mathbb{E} \left \{ \bold{w}\boldsymbol{\varphi}(\bold{w})^{\rm T} \right \}$ can be estimated from data. The estimated covariance matrix follows the Wishart distribution, based on which the number of samples can be determined to control the estimation quality.} In summary, in the context of saturated DF, the optimization problem of SMPC can be cast as follows:
\begin{equation}
\begin{split}
\min_{\bold{M},\bold{h}} &\ J(\bold{M},\bold{h}) {\rm\ defined\ in\ (\ref{eq:J1})} \\
{\rm s.t.} &\ {\rm State\ constraints\ (\ref{eq:CC})} \\
           &\ {\rm Robust\ input\ constraints\ (\ref{eq:robustinputconstraints})} \\
\end{split}
\tag{SMPC-DF}
\end{equation}

\subsection{Scenario-Based Approximation}
Due to the penetration of joint chance constraints on system states, the problem SMPC-DF is generally difficult to solve \cite{farina2016stochastic,paulson2017stochastic}. A simple yet popular technique is to replace chance constraints by hard constraints based on a finite number of realizations $\{ \bold{w}_t^{(i)} \}_{i=1}^{N_t}$, resulting in the formulation of scenario-based SMPC (SSMPC) \cite{calafiore2006scenario,calafiore2013robust}:
\begin{equation}
\begin{split}
\min_{\bold{M},\bold{h}} &\ J(\bold{M},\bold{h}) \\
\text{s.t.} &\ \bold{F}_t \left (\bold{A}x_0 + \bold{B}_u\bold{M}\boldsymbol{\varphi}(\bold{w}_t^{(i)}) + \bold{B}_u\bold{h} + \bold{B}_w\bold{w}_t^{(i)} \right ) \le \bold{f}_t,\ \forall i=1,\cdots,N_t,\ t=1,\cdots,H\ \ \ \\
&\ \sup_{|\varphi(\cdot)| \le \varphi_{\max}} \bold{G}\left [ \bold{h} + \bold{M}\boldsymbol{\varphi}(\bold{w}) \right ] \le \bold{g}.
\end{split}
\tag{SSMPC}
\end{equation}
At first glance, SSMPC is a deterministic optimization problem with a large number of constraints. Due to the randomness of sampling $N_t$ realizations, the optimal solution $\{\bold{M}^*, \bold{h}^*\}$ to SSMPC itself is a random variable. Therefore, whether $\{\bold{M}^*, \bold{h}^*\}$ satisfies the chance constraint (\ref{eq:CC}) or not is uncertain. Intuitively, if we have sufficient data samples in the scenario-based problem SSMPC, the risk that $\{\bold{M}^*, \bold{h}^*\}$ violates the chance constraint (\ref{eq:CC}) can be considerably reduced. Theoretical results have been well established in literature to demonstrate this intuition, among which a fundamental one is to choose the number of scenarios $N_t$ according to the following criterion \cite{calafiore2006scenario}:
\begin{equation}
N_t = \min \left \{ S_t \left | \sum_{j=0}^{d_t-1} \binom{N_t}{j} \epsilon^j_t (1 - \epsilon_t)^{S_t-j} \le \beta_t \right . \right \},
\label{eq:cal}
\end{equation}
where $d_t$ is the number of decision variables affecting the $t$th chance constraint and $\beta_t$ is a pre-specified confidence level. As a result, the following performance guarantee can be established in a probabilistic sense:
\begin{equation}
\mathbb{P}_{\mathcal{D}} \left \{ \mathbb{P} \left \{ \bold{F}_t \left (\bold{A}x_0 + \bold{B}_u\bold{M}^*\boldsymbol{\varphi}(\bold{w}_t) + \bold{B}_u\bold{h}^* + \bold{B}_w\bold{w}_t \right ) \le \bold{f}_t \right \} \ge 1 - \epsilon_t \right \} \ge 1 - \beta_t.
\label{eq:probguarantee}
\end{equation}
Here, the outmost probability term $\mathbb{P}_{\mathcal{D}} \left \{ \cdot \right \}$ arises from the randomness of sampling multiple scenarios, while the innermost probability $\mathbb{P}\{\cdot\}$ is associated with the disturbance sequence $\bold{w}$.

The performance guarantee (\ref{eq:probguarantee}) lies at the heart of scenario-based approximations of stochastic programs and SMPC problems. However, when DF policies are adopted, massive decision variables will be involved. In (\ref{eq:cal}), the number of decision variables induced by DF in the last stage can be calculated as:
\begin{equation}
d_H = Hn_u + n_u n_w (H-1)H/2,
\end{equation}
which leads to prohibitively massive scenarios $N_t$ that are necessitated, and hence renders the scenario-based approach impractical \cite{zhang2013stochastic}. Moreover, massive constraints in SSMPC will be encountered, which adds significant difficulty in online computations \cite{krishnamoorthy2018improving,you2018distributed}.
\section{Stochastic MPC Based on Data-Driven Robust Optimization}
In this section, we introduce the proposed SMPC approach based on data-driven robust optimization. Our basic idea is to first identify a polytope $\mathcal{W}_t(\mathcal{D})$ from available data $\mathcal{D}$ that compactly encloses the high-density region of disturbance, and then carefully \textit{calibrate} the size of $\mathcal{W}_t(\mathcal{D})$ via an independent validation dataset to ensure the probabilistic guarantee (\ref{eq:probguarantee}). Finally, a robust optimization problem is solved by enforcing constraint satisfactions for all possible realizations in the calibrated polytope.
\subsection{Uncertainty Set Learning with Support Vector Clustering}
We first introduce an efficient approach to data-based construction of $\mathcal{W}_t(\mathcal{D})$, which is proposed by \cite{shang2017data}. From a machine learning perspective, estimating the high-density region of an unknown distribution from data can be viewed as an unsupervised learning task, and SVC, which is also called one-class support vector machine, is a tailored model to this end \cite{ben2001support}. In SVC, the original data $\bold{w}$ are first mapped into a high-dimensional feature space via a nonlinear mapping $\boldsymbol{\phi}(\bold{w})$. Then, given $N$ data samples $\mathcal{D} = \{\bold{w}^{(1)},\cdots,\bold{w}^{(N)}\}$, a circle in the feature space with minimal volume is sought to include the majority of data samples. This is achieved by solving the following optimization problem:
\begin{equation}
\begin{split}
\min_{\bold{c},R,\boldsymbol{\xi}}&\ R^2 + {1 \over N\nu}\sum_{i=1}^N \xi_i\\
\mathrm{s.t.} &\ \| \boldsymbol{\phi}(\bold{w}^{(i)}) - \bold{c} \|^2 \le R^2 + \xi_i,\ i=1,\cdots,N \\
              &\ \xi_i \ge 0,\ i=1,\cdots,N
\end{split}
\label{eq:SVC_softmargin}
\end{equation}
where $\bold{c}$ and $R$ are the center and the radius of the circle, respectively. Positive slack variables $\{\xi_i\}$ are introduced to indicate data violations of the circle. If a data sample $\bold{w}^{(i)}$ lies within the circle with $\|\boldsymbol{\phi}(\bold{w}^{(i)}) - \bold{c}\|^2 \le R^2$, the associated slack variable $\xi_i$ will have zero value. Otherwise, if $\bold{w}^{(i)}$ resides outside the circle with $\|\boldsymbol{\phi}(\bold{w}^{(i)}) - \bold{c}\|^2 > R^2$, then $\xi_i > 0$. The objective in (\ref{eq:SVC_softmargin}) consists of two parts, where one minimizes the volume of the circle, and the other penalizes outliers with $\xi_i > 0$. The regularization parameter $\nu \in (0,1)$ is introduced to balance between two conflicting goals. By using Lagrange multipliers $\boldsymbol{\alpha}$, an equivalent dual problem of (\ref{eq:SVC_softmargin}) can be attained, which is essentially a quadratic program (QP) \cite{ben2001support}:
\begin{equation}
\begin{split}
\max_{\boldsymbol{\alpha}}&\ -\sum_{i=1}^N \sum_{j=1}^N \alpha_i \alpha_j K(\bold{w}^{(i)}, \bold{w}^{(j)}) + \sum_{i=1}^N \alpha_i K(\bold{w}^{(i)}, \bold{w}^{(i)}) \\
\mathrm{s.t.} &\ 0 \le \alpha_i \le 1/N\nu,\ i=1,\cdots,N \\
&\ \sum_{i=1}^N \alpha_i = 1
\end{split}
\label{eq:dualQP}
\end{equation}
Here $K(\bold{w}^{(i)}, \bold{w}^{(j)}) = \boldsymbol{\phi}(\bold{w}^{(i)})^{\rm T}\boldsymbol{\phi}(\bold{w}^{(j)})$ stands for the \textit{kernel function}, which is defined as the inner product between nonlinear mappings. Thanks to the kernel trick, one could directly specify the kernel function, without having to know the explicit form of the nonlinear mapping $\boldsymbol{\phi}(\cdot)$. After solving (\ref{eq:dualQP}), a high-density region can be established as the enclosing circle in the feature space:
\begin{equation}
\begin{split}
\mathcal{W} & = \left \{ \bold{w} \left | \|\boldsymbol{\phi}(\bold{w}) - \bold{c}\|^2 \le R^2 \right . \right \} \\
& = \left \{ \bold{w} \left | K(\bold{w}, \bold{w}) - 2 \sum_{i=1}^N \alpha_i K(\bold{w}, \bold{w}^{(i)}) + \sum_{i=1}^N \sum_{j=1}^N \alpha_i \alpha_j K(\bold{w}^{(i)}, \bold{w}^{(j)}) \le R^2 \right . \right \}.
\label{eq:uncertaintyset}
\end{split}
\end{equation}
All data points can be classified into three classes based on their different locations, namely interior points, boundary points, and outliers. According to complementary slackness, each class is characterized by specific values of associated Lagrange multipliers, as summarized in Table \ref{tab:SVCproperty}. Notice that for interior points we have $\alpha_i = 0$. Therefore, the uncertainty set (\ref{eq:uncertaintyset}) is only defined by boundary points and outliers with $\alpha_i > 0$. In this sense, $\mathcal{W}$ is ``supported" by boundary points and outliers, and this is why they are referred to as ``support vectors". For notational convenience, we define the index set of support vectors as
\begin{equation}
{\rm SV} = \left \{ i \left | \alpha_i > 0,\ 1 \le i \le N \right . \right \},
\end{equation}
and that of boundary points as
\begin{equation}
{\rm BSV} = \left \{ i \left | 0 < \alpha_i < 1/N\nu,\ 1 \le i \le N \right . \right \}.
\end{equation}

\begin{table}[t]
\caption{Classification of Data Points in Support Vector Clustering}
\centering
\begin{tabular}{c c c}
\hline
& Primal Description & Dual Description \\
\hline
Interior Points & $\|\boldsymbol{\phi}(\bold{w}^{(i)}) - \bold{c}\|^2 < R^2$ & $\alpha_i = 0$ \\
Boundary Points & $\|\boldsymbol{\phi}(\bold{w}^{(i)}) - \bold{c}\|^2 = R^2$ & $0 < \alpha_i < 1/N\nu$ \\
Outliers & $\|\boldsymbol{\phi}(\bold{w}^{(i)}) - \bold{c}\|^2 > R^2$ & $\alpha_i = 1/N\nu$ \\
\hline
\end{tabular}
\label{tab:SVCproperty}
\end{table}

As for the selection of kernel functions, the radial basis function (RBF) kernel $K(\bold{w}, \bold{v}) = \exp\{-\|\bold{w} - \bold{v}\|^2/2\sigma^2\}$ is frequently adopted in SVC literature \cite{ben2001support}. However, the RBF kernel turns out to be unsuitable for uncertainty set construction because the computational tractability of robust optimization (RO) problems will be compromised, and an analytical expression of robust counterpart problem could be unattainable. To make full use of the modeling power of SVC and preserve the tractability in RO, \cite{shang2017data} further proposed the weighted generalized intersection kernel (WGIK) that is tailored to data-driven RO:
\begin{equation}
K(\bold{w},\bold{v}) = L - \|\bold{Q}(\bold{w} - \bold{v})\|_1,
\label{eq:WGIK}
\end{equation}
where matrix $\bold{Q}$ and the scalar $L$ are kernel parameters that can be determined from data. Typically one can set $\bold{Q}$ as the sphering matrix $\bold{Q} \triangleq \boldsymbol{\Sigma}^{-\frac12}$ where $\boldsymbol{\Sigma}$ is the covariance matrix of $\bold{w}$, and value of $L$ shall be sufficiently large \cite{shang2017data}\footnote{
Different from SVC based on generic kernel functions, there is no over-fitting phenomenon for the WGIK-induced SVC model, because the uncertainty set (\ref{eq:uncertaintyset1}) will be unchanged as long as the value of $L$ is sufficiently large. Therefore, the selection of $L$ is not an issue in practice \cite{shang2017data}.}. By substituting (\ref{eq:WGIK}) into (\ref{eq:uncertaintyset}), we can arrive at the following \textit{data-driven uncertainty set}:
\begin{equation}
\mathcal{W}(\nu,\mathcal{D}) =  \left \{ \bold{w} \left | \sum_{i \in \mathrm{SV}} \alpha_i \|\bold{Q}(\bold{w} - \bold{w}^{(i)})\|_1 \le \theta \right . \right \},
\label{eq:uncertaintyset1}
\end{equation}
where
\begin{equation}
\theta = \sum_{i \in \mathrm{SV}} \alpha_i \|\bold{Q}(\bold{w}^{(i')} - \bold{w}^{(i)})\|_1,\ i' \in \mathrm{BSV}.
\end{equation}
Here we uncover the dependence of the uncertainty set $\mathcal{W}$ on both the regularization parameter $\nu$ and the dataset $\mathcal{D}$. The advantages of $\mathcal{W}(\nu,\mathcal{D})$ in RO can be demonstrated from the following aspects. First, $\mathcal{W}(\nu,\mathcal{D})$ is essentially a polytope which can be expressed as a series of linear inequalities:
\begin{equation}
\mathcal{W}(\nu,\mathcal{D}) =  \left \{ \bold{w} \left |
\begin{split}
\exists \bold{v}_i,\ & 1 \le i \le N \\
{\rm s.t.}\ \ & \sum_{i \in \mathrm{SV}} \alpha_i \bold{1}^{\rm T} \bold{v}_i \le \theta,\\
& -\bold{v}_i \le \bold{Q}(\bold{w} - \bold{w}^{(i)}) \le \bold{v}_i,\ 1 \le i \le N
\end{split}
\right . \right \},
\label{eq:uncertaintyset2}
\end{equation}
where additional auxiliary variables $\{ \bold{v}_i \}_{i=1}^N$ are introduced to eliminate 1-norm functions. Therefore, in virtue of duality theory of linear programs, either optimizing the worst-case performance or ensuring constraint satisfaction for all possible realizations in $\mathcal{W}(\nu,\mathcal{D})$ will be convenient. In other words, the tractability of RO problems with (\ref{eq:uncertaintyset2}) employed as the uncertainty set will be well preserved. Besides, since only a fraction of $\{\alpha_i\}_{i=1}^N$ will be nonzero, the induced robust counterpart problem tends to manifest a moderate complexity. Second, $\mathcal{W}(\nu,\mathcal{D})$ admits a non-parametric expression because it is defined based on SVs with $\alpha_i > 0$. Therefore, $\mathcal{W}(\nu,\mathcal{D})$ can naturally handle data asymmetry by capturing the high-density region of uncertainty from data, which is typically more compact than traditional norm-based sets and hence can help reduce the conservatism.

\subsection{Uncertainty Set Calibration}
In this subsection, we seek to endow the data-driven uncertainty set with appropriate probabilistic guarantee, which finally contributes to the performance guarantee (\ref{eq:probguarantee}) of SMPC. Motivated by \cite{hong2016approximating}, we split all available data samples into two independent datasets, namely a training dataset $\mathcal{D}_{\rm train}$ and a calibration dataset $\mathcal{D}_{\rm calib}$. After the data-driven uncertainty set $\mathcal{W}(\nu,\mathcal{D}_{\rm train})$ is developed using $\mathcal{D}_{\rm train}$, its size will be further calibrated based on $\mathcal{D}_{\rm calib}$. Note that the data-driven uncertainty set can be written as the \textit{level set} of a certain function $f(\cdot)$:
\begin{equation}
\mathcal{W}(\nu,\mathcal{D}_{\rm train}) = \left \{ \bold{w} \left | f(\bold{w}) \le \theta \right . \right \},
\label{eq:SVCset}
\end{equation}
where
\begin{equation}
f(\bold{w}) = \sum_{i \in \mathrm{SV}} \alpha_i ||\bold{Q}(\bold{w} - \bold{w}^{(i)})||_1.
\end{equation}
Hence $\theta \in \mathbb{R}$ can be considered as the ``size" parameter of the uncertainty set $\mathcal{W}(\nu,\mathcal{D}_{\rm train})$. Due to the convexity of $f(\bold{w})$, the ``volume" of $\mathcal{W}(\nu,\mathcal{D}_{\rm train})$ will become larger with the value of $\theta$ increasing. Next, the value of $\theta$ is to be carefully tuned based on an independent calibration dataset $\mathcal{D}_{\rm calib} = \{\bold{w}^{(1)}, \cdots, \bold{w}^{(N_{\rm calib})}\}$ according to the following criterion:
\begin{equation}
\tilde{\theta} = \max_{1\le i \le N_{\rm calib}} f(\bold{w}^{(i)}),
\label{eq:calibration}
\end{equation}
Finally, the calibrated data-driven set is denoted by
\begin{equation}
\tilde{\mathcal{W}}(\nu,\mathcal{D}_{\rm train}) = \left \{ \bold{w} \left | f(\bold{w}) \le \tilde{\theta} \right . \right \}.
\label{eq:calibratedset}
\end{equation}
In other words, after calibration the uncertainty set $\tilde{\mathcal{W}}(\nu,\mathcal{D}_{\rm train})$ shall enclose all data samples in the calibration dataset $\mathcal{D}_{\rm calib}$ but with a minimal volume. Now we are ready to deal with the probabilistic guarantee of the data-driven uncertainty set. We first restate a useful lemma in literature.

\noindent \textbf{Lemma 1} \textit{(Theorem 8.1 in \cite{tempo2012randomized})}. Let $\{ Y_1,\cdots,Y_N \}$ be i.i.d. samples from a random variable $Y \in \mathbb{R}$. For any $\epsilon \in (0,1)$ and $\beta \in (0,1)$, if
\begin{equation}
N \ge \frac{\log \beta}{\log (1 - \epsilon)},
\end{equation}
then, with probability greater than $1 - \beta$ we have:
\begin{equation}
\mathbb{P} \left \{ Y \le \max_{i=1,\cdots,N} Y_i \right \} \ge 1 - \epsilon.
\end{equation}
That is,
\begin{equation}
\mathbb{P}_{\mathcal{D}} \left \{ \mathbb{P} \left \{ Y \le \max_{i=1,\cdots,N} Y_i \right \} \ge 1 - \epsilon \right \} \ge 1 - \beta.
\end{equation}

Next, the performance guarantee can be established in terms of the following theorem, as a major contribution of this work.

\noindent \textbf{Theorem 1:} If the parameter $\tilde{\theta}$ in (\ref{eq:calibratedset}) is determined according to (\ref{eq:calibration}) and the number of independent data samples $N_{\rm calib}$ satisfy
\begin{equation}
N_{\rm calib} = \lceil \log \beta / \log (1-\epsilon) \rceil,
\label{eq:newresult}
\end{equation}
then the calibrated SVC-based uncertainty set admits the following property:
\begin{equation}
\mathbb{P}_{\mathcal{D}} \left \{ \mathbb{P} \left \{ \bold{w} \in \tilde{\mathcal{W}}(\nu,\mathcal{D}_{\rm train}) \right \} \ge 1 - \epsilon \right \} \ge 1 - \beta.
\label{eq:probgrt}
\end{equation}
\textit{Proof:} Here we regard $f(\bold{w})$ as a random variable. In the spirit of Lemma 1, the maximum of $\{ f(\bold{w}^{(i)}) \}_{i=1}^{N_{\rm calib}}$ is a valid $1-\beta$ confidence upper bound for the $(1-\epsilon)$-quantile of the random variable $f(\bold{w})$, which is regarded as a random variable. That is,
\begin{equation}
\mathbb{P}_{\mathcal{D}} \left \{ \mathbb{P} \left \{ f(\bold{w}) \le \tilde{\theta} = \max_{1\le i \le N_{\rm calib}} f(\bold{w}^{(i)}) \right \} \ge 1 - \epsilon \right \} \ge 1 - \beta,
\end{equation}
which turns out to be equal to (\ref{eq:probgrt}). This completes the proof. \qed

\vspace{12pt}

\noindent \textbf{Remark 1}. The regularization parameter $\nu$ in SVC can be selected as $\nu = \epsilon$. The rationality lies in that $1 - \nu$ is an lower bound on the fraction of data coverage on training samples \cite{shang2017data}.

\subsection{A Data-Driven Robust Optimization Scheme}
For constraints in the $t$th stage in SSMPC, we first adopt the lifting technique \cite{georghiou2015generalized} to define the lifted uncertainty:
\begin{equation}
\hat{\bold{w}}_t = \left [
\begin{array}{c}
\boldsymbol{\varphi}(\bold{w}_t) \\
\bold{w}_t \\
\end{array}
\right ] \in \mathbb{R}^{2tn_w},\ t = 1,\cdots,H.
\end{equation}
Then constraints induced by $N_t$ scenarios $\{\bold{w}_t^{(i)}\}_{i=1}^{N_t}$ can be recast as being linear in the lifted uncertainty:
\begin{equation}
\bold{F}_t \left (\bold{A}x_0 + [\bold{B}_u\bold{M}\ \bold{B}_w]\hat{\bold{w}}_t^{(i)} + \bold{B}_u\bold{h} \right ) \le \bold{f}_t,\ \forall i=1,\cdots,N_t
\end{equation}
where the lifted sample
\begin{equation}
\hat{\bold{w}}_t^{(i)} = \left [
\begin{array}{c}
\boldsymbol{\varphi}(\bold{w}_t^{(i)}) \\
\bold{w}_t^{(i)} \\
\end{array}
\right ]
\end{equation}
is formulated based on the $i$th scenario $\bold{w}_t^{(i)},\ \forall i=1,\cdots,N_t$. We use a training dataset $\hat{\mathcal{D}}_{\rm train}^{(t)} = \{ \hat{\bold{w}}_t^{(i)} \}_{i=1}^{N^{(t)}_{\rm train}}$, which includes $N_{\rm train}^{(t)}$ samples of lifted uncertainty, to establish the SVC-based uncertainty set, and use $\hat{\mathcal{D}}_{\rm calib}^{(t)}$, which includes $N_{\rm calib}^{(t)}$ samples of lifted uncertainty, for calibration with $N_{\rm calib}^{(t)} \ge \log \beta_t / \log (1-\epsilon_t)$, thereby finally yielding the calibrated dataset $\tilde{\mathcal{W}}_{t} ( \nu, \hat{\mathcal{D}}_{\rm train}^{(t)} )$.
Finally, we solve the following data-driven robust optimization (DDRO) problem:
\begin{equation}
\begin{split}
\min_{\bold{M},\bold{h}} &\ J(\bold{M},\bold{h}) \\
{\rm s.t.} &\ \bold{F}_t(\bold{A}x_0 + [\bold{B}_u\bold{M}\ \bold{B}_w]\hat{\bold{w}}_t + \bold{B}_u\bold{h} ) \le \bold{f}_t,\ \forall \hat{\bold{w}}_t \in \tilde{\mathcal{W}}_{t} ( \nu, \hat{\mathcal{D}}_{\rm train}^{(t)} ),\ t=1,\cdots,H \ \ \ \ \ \ \ \ \\
&\ \sup_{|\varphi(\cdot)| \le \varphi_{\max}} \bold{G}\left [ \bold{h} + \bold{M}\boldsymbol{\varphi}(\bold{w}) \right ] \le \bold{g}
\end{split}
\tag{DDRO}
\end{equation}
and the resulted predictive controller is referred to as the data-driven robust MPC (DRMPC) hereafter. Since $\tilde{\mathcal{W}}_{t} ( \nu, \hat{\mathcal{D}}_{\rm train}^{(t)} )$ can be expressed as a series of linear inequalities, and $\sup_{|\varphi(\cdot)| \le \varphi_{\max}} \bold{G}\left [ \bold{h} + \bold{M}\boldsymbol{\varphi}(\bold{w}) \right ] \le \bold{g}$ can be regarded as robust linear constraints with polytopic uncertainty sets \cite{paulson2017stochastic}, DDRO is essentially a classic robust optimization problem and hence can be desirably transformed into a tractable convex program \cite{bertsimas2004price}. Moreover, we can readily arrive at the following theorem.

\noindent \textbf{Theorem 2}. The optimal solution $\{\bold{M}^*,\bold{h}^*\}$ to DDRO satisfies the probabilistic guarantee (\ref{eq:probguarantee}).

\noindent \textit{Proof:} With a feasible solution $\{ \bold{M}, \bold{h} \}$ satisfying constraints in DDRO, the event

\noindent $\mathbb{P} \left \{ \hat{\bold{w}}_t \in \tilde{\mathcal{W}}_t(\nu,\hat{\mathcal{D}}^{(t)}_{\rm train}) \right \} \ge 1 - \epsilon_t$ implies $\mathbb{P} \left \{ \bold{F}_t(\bold{A}x_0 + [\bold{B}_u\bold{M}\ \bold{B}_w]\hat{\bold{w}}_t + \bold{B}_u\bold{h} ) \le \bold{f}_t \right \} \ge 1 - \epsilon_t$. Then we have
\begin{equation}
\begin{split}
&\ \mathbb{P}_{\mathcal{D}} \left \{ \mathbb{P} \left \{ \bold{F}_t(\bold{A}x_0 + [\bold{B}_u\bold{M}\ \bold{B}_w]\hat{\bold{w}}_t + \bold{B}_u\bold{h} ) \le \bold{f}_t \right \} \ge 1 - \epsilon_t \right \} \\
\ge &\ \mathbb{P}_{\mathcal{D}} \left \{ \mathbb{P} \left \{ \hat{\bold{w}}_t \in \tilde{\mathcal{W}}_t(\nu,\hat{\mathcal{D}}^{(t)}_{\rm train}) \right \} \ge 1 - \epsilon_t \right \} \\
\ge &\ 1 - \beta_t.
\end{split}
\label{eq:47}
\end{equation}
Because the optimal solution $\{\bold{M}^*,\bold{h}^*\}$ is feasible for DDRO, (\ref{eq:47}) directly yields the probabilistic guarantee (\ref{eq:probguarantee}). This completes the proof. \qed

\vspace{12pt}

The above theorem indicates that, by addressing the DDRO problem we \textit{approximately} solve the problem of SMPC with appropriate performance guarantee preserved. Besides, as a systematic data-driven approach, the proposed DRMPC possesses some inherent advantages, as discussed below.

\noindent \textbf{Remark 2}. The proposed DRMPC approach applies not only to the cases with the SVC-based uncertainty set employed. It can be easily shown that the calibration procedure still takes effect with other machine learning techniques used to construct the data-driven uncertainty set, provided that the uncertainty set can be expressed as the level set of a function $\mathcal{W}(\mathcal{D}) = \{ \bold{w} | f(\bold{w}) \le \theta \}$.

\noindent \textbf{Remark 3}. If $\epsilon_t \equiv \epsilon$ and $\beta_t \equiv \beta$ for all stages, it suffices to construct only one uncertainty set $\tilde{\mathcal{W}}_H(\nu,\hat{\mathcal{D}}^{(H)}_{\rm train})$ for the last stage and then calibrate it with $\lceil \log \beta / \log (1-\epsilon) \rceil$ data samples. The uncertainty sets $\tilde{\mathcal{W}}_t(\nu,\hat{\mathcal{D}}^{(t)}_{\rm train}),\ t=1,\cdots,H-1$ for earlier stages can be easily obtained by projecting $\tilde{\mathcal{W}}_H(\nu,\hat{\mathcal{D}}^{(H)}_{\rm train})$ onto a portion of dimensions, which saves potential computational costs but still preserves probabilistic guarantees.

\noindent \textbf{Remark 4}. Although data-driven uncertainty sets in DDRO are built based on scenarios, most scenarios are characterized by zero Lagrange multipliers $\{\alpha_i = 0\}$, and only support vectors have non-zero Lagrange multipliers $\{ \alpha_i > 0 \}$, whose portion is approximately controlled by the regularization parameter $\nu$. Therefore, the formulation of the data-driven uncertainty set (\ref{eq:uncertaintyset2}) will not be over-complicated \cite{shang2017data}.

\noindent \textbf{Remark 5}. It is worth mentioning that, the proposed method cannot be applied to the case of multiplicative uncertainty, while classic SSMPC is able to handle multiplicative uncertainty. For predictive control of general nonlinear systems with both multiplicative and additive uncertainty, one can utilize the polynomial chaos expansion-based and moment-based formulations of SMPC \cite{paulson2017input,paulson2017efficient}.

\subsection{Discussions}
The proposed approach to SMPC shares some similarities with the method in \cite{zhang2013stochastic} combining scenario optimization with robust optimization, where a hyper-rectangle capturing all scenarios $\mathcal{W}_t(\mathcal{D}) = \left \{ \bold{w}_t | \boldsymbol{\gamma}_{\min}^* \le \bold{w}_t \le \boldsymbol{\gamma}_{\max}^* \right \}$ is first identified by solving the following scenario-based optimization problem:
\begin{equation}
\begin{split}
\min_{\boldsymbol{\gamma}_{\min}, \boldsymbol{\gamma}_{\max}} & \ \|\boldsymbol{\gamma}_{\max} - \boldsymbol{\gamma}_{\min}\|_1 \\
{\rm s.t.}\ \ \ & \ \ \boldsymbol{\gamma}_{\min} \le \bold{w}_t^{(i)} \le \boldsymbol{\gamma}_{\max},\ \forall i = 1,\dots,N_t
\end{split}
\label{eq:aux}
\end{equation}
Notice that there are only $2tn_w$ decision variables in (\ref{eq:aux}); therefore, if $N_t$ satisfies (\ref{eq:cal}) with $d_t = 2tn_w$, the optimal solution to (\ref{eq:aux}) admits the following probabilistic guarantee:
\begin{equation}
\mathbb{P}_{\mathcal{D}} \left \{ \mathbb{P} \left \{ \bold{w}_t \in \mathcal{W}_t(\mathcal{D}) \right \} \right \} = \mathbb{P}_{\mathcal{D}} \left \{ \mathbb{P} \left \{ \boldsymbol{\gamma}_{\min}^* \le \bold{w}_t \le \boldsymbol{\gamma}_{\max}^* \right \} \ge 1 - \epsilon_t \right \} \ge 1 - \beta_t.
\label{eq24}
\end{equation}
The crux to a reduced sampling complexity lies in that, there are much fewer decision variables in (\ref{eq:aux}) than in generic SSMPC, that is, $2tn_w < tn_u + n_u n_w (t-1)t/2$ when $t$ is sufficiently large. Finally, chance constraints (\ref{eq:CC}) are replaced by the following robust constraints:
\begin{equation}
\bold{F}_t(\bold{A}x_0 + \bold{B}_u\bold{M}\boldsymbol{\varphi}\left ( \bold{w}_t \right ) + \bold{B}_u\bold{h} + \bold{B}_w\bold{w}_t) \le \bold{f}_t,\ \forall \bold{w}_t \in \mathcal{W}_t,\ t=1,\cdots,H
\label{eq:RCC}
\end{equation}
which leads to an RMPC problem.

Despite the established probabilistic guarantee, the strategy proposed in \cite{zhang2013stochastic} is still prone to some inherent deficiencies. On one hand, the sampling size required by \cite{zhang2013stochastic} increases with stage $t$. In this regard, we make the following remarks to highlight the advantages of the proposed DRMPC approach.

\noindent \textbf{Remark 6}. The calibration sample size required by DRMPC depends only on the violation probability $\epsilon$ and the confidence level $\beta$, and it is unrelated to the dimension of uncertainty and the horizon length $H$. In contrast, the sample size required by \cite{zhang2013stochastic} is dependent on the dimension of uncertainty, which tends to be high when a long prediction horizon is adopted.

\noindent \textbf{Remark 7}. The total sample size required by DRMPC is $N_{\rm train} + N_{\rm calib}$. Notice that $N_{\rm train}$ training samples only affect the ``initial shape" of the data-driven uncertainty set, which can be characterized by the function $f(\bold{w})$, and the probabilistic guarantee is ensured whenever $N_{\rm calib}$ fulfills (\ref{eq:newresult}). In principle, one shall not use very few training samples in order to reduce the overall sampling complexity, although the probabilistic guarantee still holds with sufficient calibration data. This is because too few training samples may lead to under-fitting of the high-density region of uncertain parameter space, which further induces over-conservatism. Therefore, the advantages of the proposed DRMPC approach shall be understood as that, by choosing a moderate number of training samples $N_{\rm train}$, it is probable in practice to learn a compact envelope of uncertainty in the training stage, while still achieving a reduced overall sampling complexity compared with \cite{calafiore2006scenario,zhang2013stochastic}.

\subsubsection*{A Motivating Example}

\begin{figure}[h]
\centering
\subfigure[]{
\includegraphics[width = 0.4\textwidth]{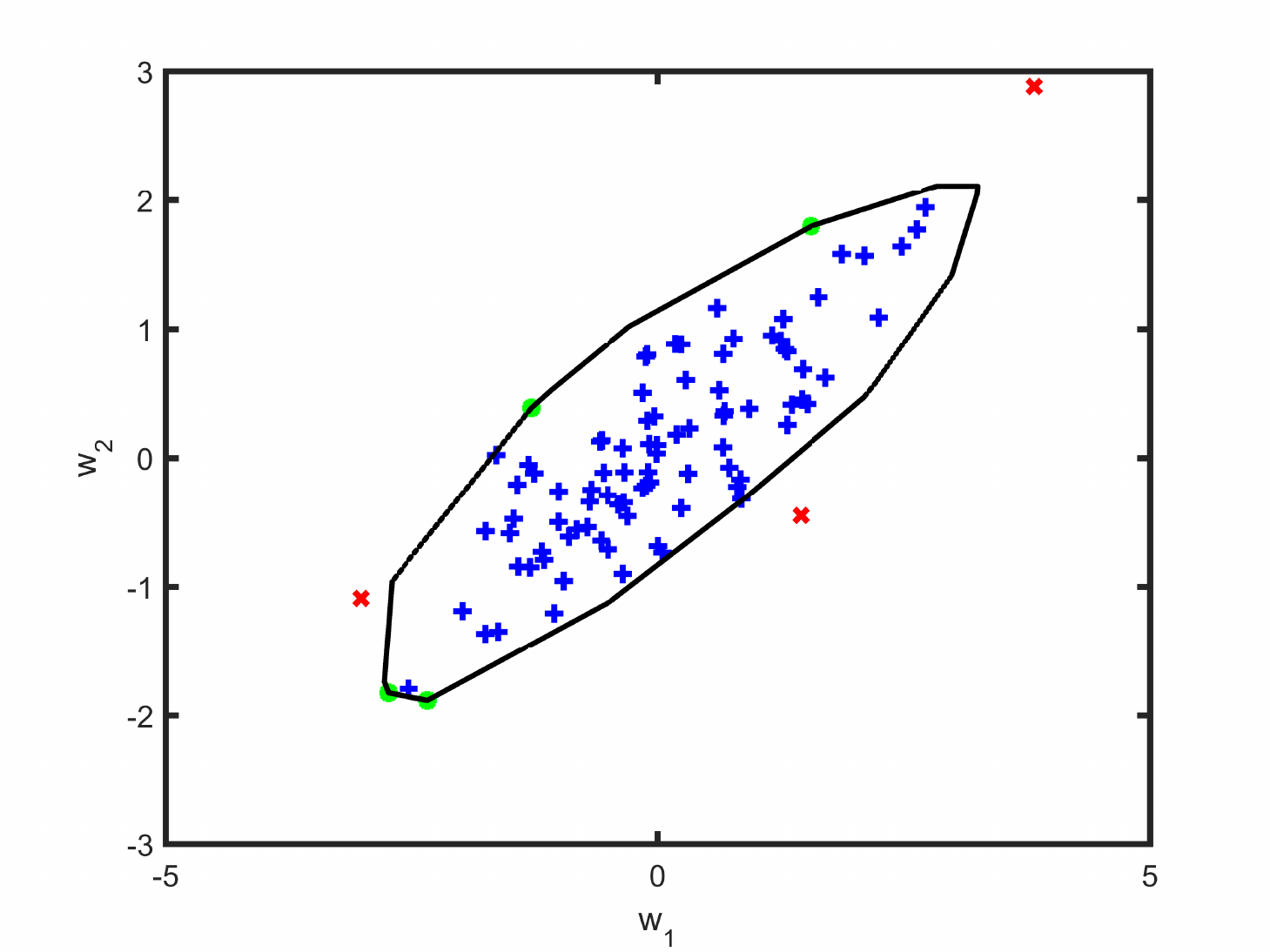}}
\subfigure[]{
\includegraphics[width = 0.4\textwidth]{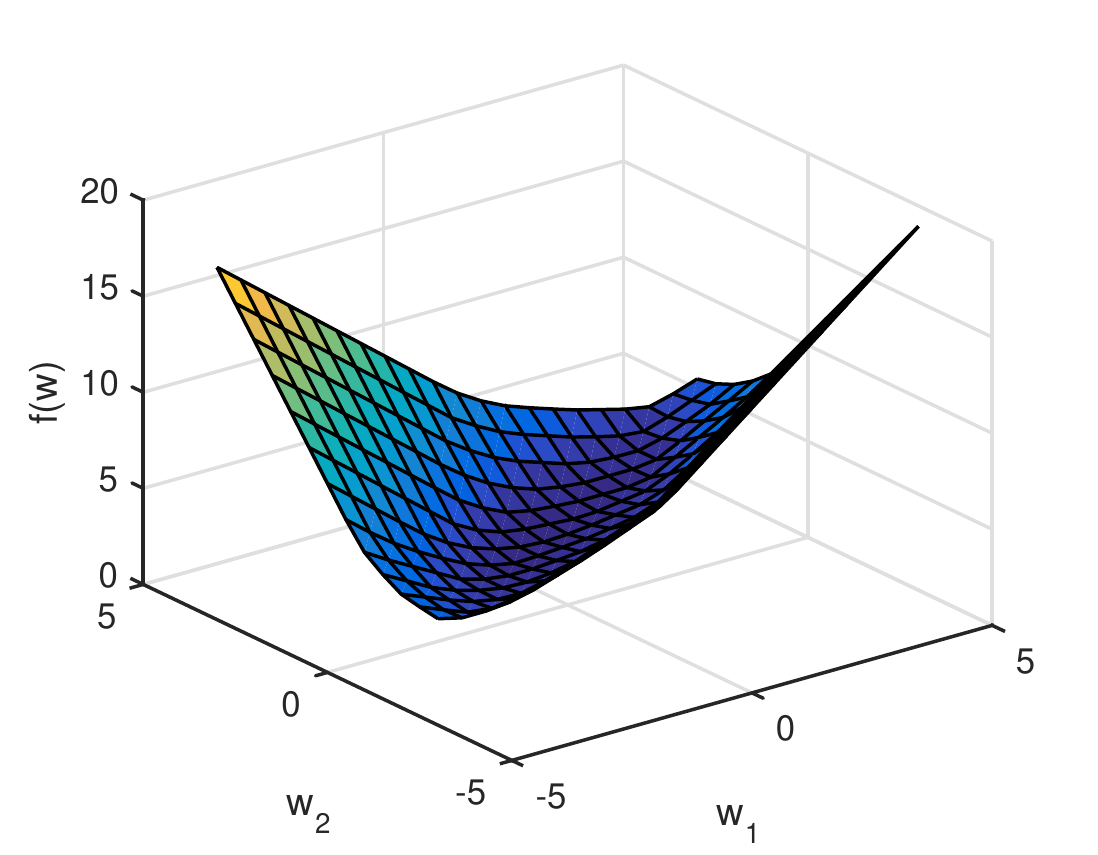}}
\subfigure[]{
\includegraphics[width = 0.4\textwidth]{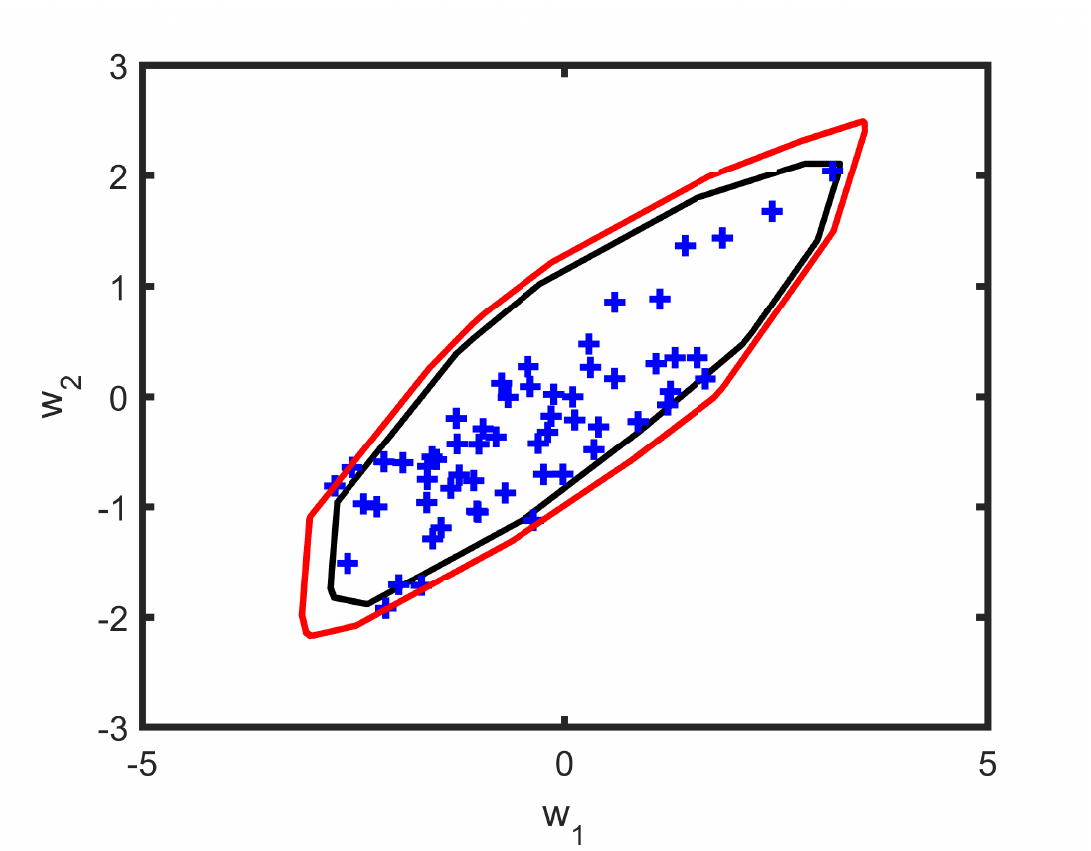}}
\subfigure[]{
\includegraphics[width = 0.4\textwidth]{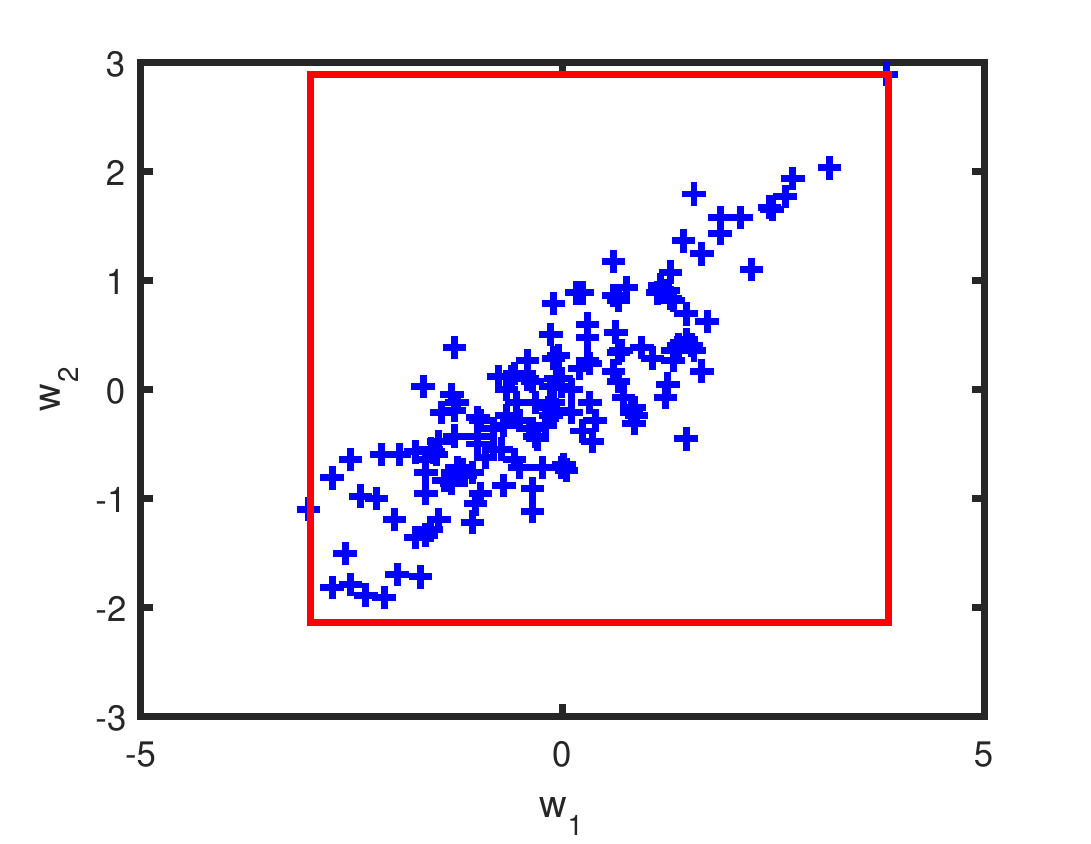}}
\caption{A 2-D example of the data-driven uncertainty set construction and calibration procedure. (a) Training data and the uncertainty set. Interior points, boundary points, and outliers are marked as blue, green, and red, respectively. (b) The geometric shape of $f(\bold{w})$. (c) Uncertainty set calibration based on a calibration dataset. The original and the calibrated uncertainty sets are marked as black and red, respectively. (d) Hyper-rectangle uncertainty set.}
\label{fig:2d_example}
\end{figure}

Next, a simple yet illustrative example is provided to further explain the merit of the SVC-based uncertainty set. Assume that the uncertainty $\bold{w} \in \mathbb{R}^2$ follows a bivariate Gaussian distribution. We collect $N_{\rm train} = 94$ scenarios to construct the training dataset $\mathcal{D}_{\rm train}$, and then build the SVC-based uncertainty set with $\nu = 0.05$. The results are shown in Figure \ref{fig:2d_example}(a), where it can be observed that a polytope is capable of capturing the geometric shape of the distribution and including most data samples compactly. Figure \ref{fig:2d_example}(b) visualizes the function $f(\bold{w})$, which is a convex piesewise linear function. In the calibration step, we use the parameter setting $\epsilon = \beta = 0.05$. This results in $N_{\rm calib} = 59$ according to (\ref{eq:newresult}), and the value of $\theta$ is calibrated based on an independent calibration dataset $\mathcal{D}_{\rm calib}$, as shown in Figure \ref{fig:2d_example}(c). Here 59 samples for calibrations have been plotted. Before calibration $\theta = 4.1146$, and after calibration the value of $\tilde{\theta}$ becomes 4.2714. Therefore, the new uncertainty set (marked with red edges) becomes larger and encloses all calibration data successfully.

For comparison reasons, the hyper-rectangle uncertainty set in \cite{zhang2013stochastic} is also reported in Figure \ref{fig:2d_example}(d). To ensure the probabilistic guarantee with $\epsilon = \beta = 0.05$, $N = 153$ samples are needed according to (\ref{eq:cal}) with $d_t = 2n_w = 4$. In this way, the same amount of data points are used in two different approaches for uncertainty set constructions (for the SVC-based set $N_{\rm train} + N_{\rm calib} = 153$), which provide identical probabilistic guarantees. It can be seen from Figure \ref{fig:2d_example}(d) that the hyper-rectangle uncertainty set gives rise to unnecessary coverage, which probably leads to over-conservative solutions of RMPC problems. In contrast, the SVC-based uncertainty set enables a fairly compact coverage. In reality, the distribution $\mathbb{P}\{\cdot\}$ of disturbance $\bold{w}_t$ probably embodies some intrinsic structure, such as cross-correlation and asymmetry, which cannot be appropriately handled by hyper-rectangular set but can be addressed by SVC-based uncertainty set. For instance, in optimal control and operations of energy systems \cite{velloso2018scenario} and agricultural systems \cite{delgoda2016irrigation}, one typically desires to hedge against uncertainty of future weather forecast, which embodies temporally correlations but allows for historical data accumulation. Therefore, by accommodating intrinsic structures within the uncertainty data, the SVC-based uncertainty sets can help reducing the conservatism of control decisions.

\subsection{Feasibility and Stability}
In this subsection, we delve into the issue of feasibility and stability of the proposed control scheme. We first define the set of initial states $x$ for which a feasible saturated DF policy exists as:
\begin{equation}
X^{\rm DF}_H \triangleq \left \{ x_0 \in \mathbb{R}^{n_x} \left | {\rm\ DDRO\ with\ initial\ states\ } x_0 {\rm \ is\ feasible} \right . \right \}.
\end{equation}
Due to the boundedness of system inputs $\{ u_t \}$, $X^{\rm DF}_H$ is obviously bounded in the presence of unbounded disturbances. Therefore, although constraint satisfactions are enforced on the data-driven uncertainty set in the optimization problem DDRO, it is often the case that $x_0 \notin X^{\rm DF}_H$, thereby compromising feasibility of DDRO. It implies that a robustly invariant set, which has been created to ensure recursive feasibility in robust RHC literature \cite{mayne2000constrained,goulart2006optimization}, no longer exists. Therefore, we solve a \textit{backup control problem} that always remains feasible when system states fail to reside in $X^{\rm DF}_H$. This can be made, for instance, by softening the state constraints and absorbing them in the objective by means of an exact penalty function \cite{calafiore2013robust}. Here we do not focus on a particular approach since a number of effective strategies have been developed in literature. In this way, even if DDRO is infeasible a backup control decision $\bar{\mu}(x_t)$ can still be reliably resorted to. Table \ref{tab:procedure} summarizes the implementation procedure of the RHC scheme with an additional backup controller.

\begin{table}[htbp]
\caption{Implementation Procedure with A Backup Controller to Ensure Feasibility}
\centering
\begin{tabular}{l}
\hline
At every time stage $t$, carry out the following steps: \\
\textit{Step 1}. Measure the current system states $x_t$. \\
\textit{Step 2}. Solve the problem DRMPC with initial states $x_t$. If it is feasible, obtain the \\
\ \ \ \ \ \ \ \ \ \ \ RHC control law $u_t = h_0^*$, and go to Step 4; Otherwise, go to Step 3.\\
\textit{Step 3}. Solve the backup control problem, and get the RHC control law $u_t = \bar{\mu}(x_t)$. \\
\textit{Step 4}. Implement the control action $u_t$ onto the plant. \\
\hline
\end{tabular}
\label{tab:procedure}
\end{table}

The stability issue of stochastic linear systems under RHC laws has been investigated by \cite{chatterjee2011stochastic,cherukuri2011stochastic}, which can be briefly described as $\sup_{t \in \mathbb{N}_0} \mathbb{E} \{ \|x_t\|^2 \} < \infty$; however, their results bear apparent limitations in practice because of the following reasons. In \cite{chatterjee2011stochastic}, no state constraints are involved, and the resulting optimization problem is globally feasible, whereas the work by \cite{cherukuri2011stochastic} does consider state constraints but assumes the disturbances to be bounded. In this work, we consider a more comprehensive setting where the RHC law is adopted, state constraints are present, and the disturbance is unbounded and non-i.i.d.

The stability of the stochastic linear system (\ref{eq:statespace}) substantially depends on the property of system matrix $A$, where explicit conclusions can be drawn in three distinct categories. If $A$ is Schur-stable, with all eigenvalues inside the unit circle, the mean-square boundedness $\sup_{t \in \mathbb{N}_0} \mathbb{E} \{ \|x_t\|^2 \} < \infty$ holds with bounded inputs and the bounded covariance $\mathbb{E} \{ \bold{w}_t \bold{w}_t^{\rm T} \}$ of disturbances \cite{paulson2017stochastic}. If $A$ has at least one eigenvalue outside the unit circle, it has been shown that mean-square boundedness cannot be attained with a bounded control law. When $A$ is Lyapunov stable (eigenvalues on the unit circle with equal geometric and algebraic multiplicities), it suffices to consider the situation when $A$ is orthogonal ($A^{\rm T}A = I$) \cite{chatterjee2011stochastic}. To tackle this case, we state that, the stability issue has an intimate relationship to the recursive feasibility issue, especially the performance of the backup control policy $\bar{\mu}(x_t)$.

First, some additional assumptions are made as follows.

\noindent \textbf{Assumption 1}. The stochastic process $\{ w_t \}_{t \in \mathbb{N}_0}$ satisfies
\begin{equation}
\sup_{t \in \mathbb{N}_0} \mathbb{E} \{ \| w_t \|^4 \} = C_4 < \infty.
\end{equation}

\noindent \textbf{Assumption 2}. For the optimal control policy $\bar{\mu}(x_t)$ attained by the backup control problem in Table \ref{tab:procedure}, there exists a constant $r > \| B_w \|_2 \cdot \sup_{t \in \mathbb{N}_0} \mathbb{E} \{ \| w_t \| \}$ such that
\begin{equation}
\|Ax_t + B_u \bar{\mu}(x_t)\| - \|x_t\| \le -r,\ \forall x_t \notin X^{\rm DF}_H.
\label{eq:assumption2}
\end{equation}

Assumption 1 requires $\{ w_t \}_{t \in \mathbb{N}_0}$ to have a finite fourth moment, which is not restrictive. Assumption 2 speaks about the performance of the backup controller in steering the system states towards origin when DDRO is infeasible. An intuitive interpretation is that, a minimal improvement, which is measured by $r$, is guaranteed to be made upon the \textit{nominal system states} at the next sampling instant $t+1$, and the improvement must exceed the largest expectation of $\| B_w \|_2\cdot\|w_t\|$. In a nutshell, Assumption 2 requires that the effect of the backup control in steering system states towards the origin, must dominate the effect of additive disturbances. In this sense, Assumption 2 is not restrictive in practice. Next, we introduce a fundamental result in stochastic processes that will be useful in the sequel.

\noindent \textbf{Lemma 2} \textit{\cite{pemantle1999moment}}. Let $\{ \eta_t \}_{t \in \mathbb{N}_0}$ be a sequence of nonnegative random variables on some probability space $\left ( \Delta, \mathcal{F}, \mathbb{P} \right )$, and let $\{ \mathcal{F}_t \}_{t \in \mathbb{N}_0}$ be any filtration to which the stochastic process $\{ \eta_t \}_{t \in \mathbb{N}_0}$ is adapted. Suppose that there exist constants $b > 0$, and $Z, M < \infty$, such that $\eta_0 < Z$, and for all $t$,
\begin{gather}
\mathbb{E} \left \{ \eta_{t + 1} - \eta_t | \mathcal{F}_t \right \} \le -b\ \mathrm{on\ the\ event} \left \{ \eta_t > Z \right \},\ \mathrm{and} \label{eq:lemma21} \\
\mathbb{E} \left \{ |\eta_{t + 1} - \eta_t|^4 | \eta_0, \cdots, \eta_t \right \} \le M. \label{eq:lemma22}
\end{gather}
Then there exists a constant $\gamma = \gamma(b,Z,M) > 0$ such that $\sup_{t \in \mathbb{N}_0} \mathbb{E} \{ \eta_t^2 \} <\gamma$.

Details of the proof can be found in \cite{pemantle1999moment} and are omitted here. Now we are ready to present the main stability result.

\noindent \textbf{Theorem 3}. Suppose that admissible control inputs are bounded $\sup_{u \in F_u} \|u\| = r_u < \infty$ and $A$ is orthogonal, and Assumptions 1 and 2 are satisfied, then for all initial states $x_0 \in \Pi_H(x)$ there always exists a constant $\gamma > 0$ such that the closed-loop system admits the mean-square boundedness:
\begin{equation}
\sup_{t \in \mathbb{N}_0} \mathbb{E} \{ \|x_t\|^2 \} = \gamma < \infty.
\label{eq:stability}
\end{equation}
\textit{Proof:} The condition (\ref{eq:lemma21}) is to be verified at first. A non-negative stochastic process is defined as $\{ \eta_t = \|x_t\| \}_{t \in \mathbb{N}_0}$. Since $X^{\rm DF}_H$ is bounded, we set $Z = \sup_{x \in X^{\rm DF}_H} \|x\| < \infty$. On the event $\eta_t \ge Z$, $x_t \notin X^{\rm DF}_H$ is infeasible for DDRO, and in this case the backup control problem is to be solved. In this case, we have
\begin{equation}
\begin{split}
&\ \mathbb{E} \left \{ \eta_{t + 1} - \eta_t | \mathcal{F}_t \right \} \\
= &\ \mathbb{E} \left \{ \|Ax_t + B_u \bar{\mu}(x_t) + B_w w_t\| - \|x_t\| \left | \mathcal{F}_t \right . \right \} \\
\le &\ \mathbb{E} \left \{ \|Ax_t + B_u \bar{\mu}(x_t)\| + \|B_w w_t\| - \|x_t\| \left | \mathcal{F}_t \right . \right \} \\
\le &\ -r + \mathbb{E} \left \{ \|B_w w_t\| \left | \mathcal{F}_t \right . \right \}\ \ \  ({\rm Assumption\ 2})\\
\le &\ -r + \|B_w\|_2 \cdot \sup_{\|w_t\|} \mathbb{E} \left \{ \|w_t\| \right \}
\end{split}
\end{equation}
According to Assumption 2, we can define a positive scalar $b = r - \|B_w\|_2 \cdot \sup_{\|w_t\|} \mathbb{E} \left \{ \|w_t\| \right \} > 0$, which fulfills (\ref{eq:lemma21}). Next we deal with (\ref{eq:lemma22}) by defining new system states $y_t = \left ( A^{\rm T} \right )^t x_t$. Then we have
\begin{equation}
\begin{split}
y_{t + 1} = & \left ( A^{\rm T} \right )^{t+1} x_{t + 1} \\
= &\ \left ( A^{\rm T} \right )^{t} x_{t} + \left ( A^{\rm T} \right )^{t+1} \left [ B_u \bar{\mu}(x_t) + B_w w_t \right ] \\
= &\ y_t + \left ( A^{\rm T} \right )^{t+1} \left [ B_u \bar{\mu}(x_t) + B_w w_t \right ]
\end{split}
\end{equation}
Meanwhile, because $\|y_t\|^2 = y_t^{\rm T}y_t = x_t^{\rm T}A^{t} (A^{\rm T})^{t}x_t = x_t^{\rm T}x_t = \|x_t\|^2$, we have $\|y_t\| = \|x_t\|$. According to the triangle inequality $-\|y_{t+1} - y_t\| \le \|y_{t+1}\| - \|y_{t}\| \le \|y_{t+1} - y_t\|$, it holds that $| \|y_{t+1}\| - \|y_{t}\| |^4 \le \|y_{t+1} - y_t\|^4$, which indicates that:
\begin{equation}
\begin{split}
& |\eta_{t+1} - \eta_t|^4 \\
= & \lvert \|y_{t+1}\| - ||y_{t}\| \rvert ^4 \\
\le & \|y_{t+1} - y_t\|^4 \\
= & \|\left ( A^{\rm T} \right )^{t+1} \left [ B_u \bar{\mu}(x_t) + B_w w_t \right ]\|^4 \\
= & \| B_u \bar{\mu}(x_t) + B_w w_t \|^4 \\
\le & \left ( \|B_u\|_2 \cdot \|\bar{\mu}(x_t)\| + \|B_w\|_2 \cdot \|w_t\| \right )^4 \\
\le & \left ( \|B_u\|_2 \cdot r_u + \|B_w\|_2 \cdot \|w_t\| \right )^4
\end{split}
\end{equation}
Therefore, $\mathbb{E} \left \{ |\eta_{t + 1} - \eta_t|^4 \left | \eta_0, \cdots, \eta_t \right . \right \} \le \mathbb{E} \left \{ \left ( \|B_u\|_2 \cdot r_u + \|B_w\|_2 \cdot \|w_t\| \right )^4 \right \}$. Because of the boundedness of the fourth moment of $\|w_t\|$, it can be easily deduced that there exists a constant $M = M(\|B_u\|_2, r_u, \|B_w\|_2, C_4)>0$ such that
\begin{equation}
\mathbb{E} \left \{ \left ( \|B_u\|_2 \cdot r_u + \|B_w\|_2 \cdot \|w_t\| \right )^4 \right \} \le M,
\end{equation}
which yields the condition (\ref{eq:lemma22}). Now all constants $\{b,Z,M\}$ are well-defined. Note that $x_0 \in \Pi_H(x)$ implies $\eta_0 = \|x_0\| < Z$. Therefore, in the light of Lemma 2, there exists a constant $\gamma = \gamma(b,Z,M) > 0$ such that
\begin{equation}
\sup_{t \in \mathbb{N}_0} \mathbb{E} \{ \|x_t\|^2 \}  = \sup_{t \in \mathbb{N}_0} \mathbb{E} \{ \eta_t^2 \} <\gamma.
\end{equation}
This completes the proof. \qed

\section{Case Study on A Two-Mass-Spring System}
\begin{figure}[h]
\centering
\includegraphics[width=0.45\textwidth]{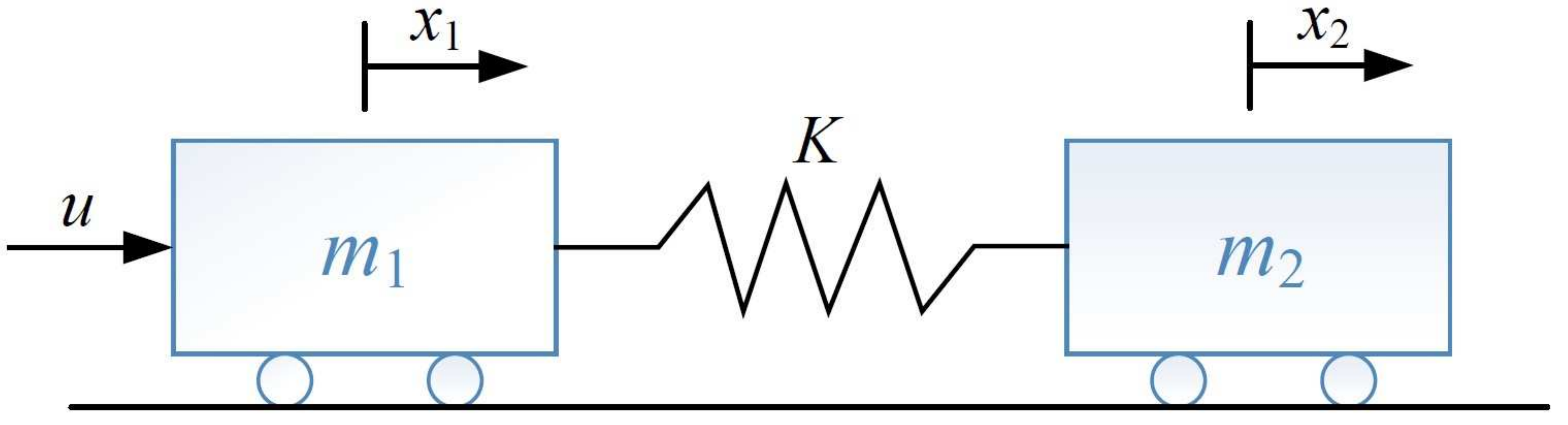}\\
\caption{A two-mass-spring system \cite{kothare1996robust}.}
\label{fig:mass_spring}
\end{figure}

In this section, we adopt a benchmark two-mass-spring system from \cite{kothare1996robust} and carry out a comparative study between generic SMPC methods and the proposed method. The system consists of two masses connected by a spring, as shown in Figure \ref{fig:mass_spring}. The dynamics of the system can be approximated by (\ref{eq:statespace}) with the sampling time of 0.1s, where system matrices are given by:
\begin{equation}
A = \left [
\begin{array}{c c c c}
1 & 0 & 0.1 & 0 \\
0 & 1 & 0 & 0.1 \\
-\frac{K}{m_1} & 0.1\frac{K}{m_1} & 1 & 0 \\
\frac{K}{m_2} & -0.1\frac{K}{m_2} & 0 & 1 \\
\end{array}\right ],\
B_u = \left [
\begin{array}{c}
0 \\
0 \\
0.1\frac{1}{m_1} \\
0 \\
\end{array}
\right ],\ B_w = \left [
\begin{array}{c}
1.0 \\
0.5 \\
0.3 \\
0.4 \\
\end{array}
\right ].
\end{equation}
Here, model states $x_t$ have four dimensions, where the first two states $\{x_t(1), x_t(2)\}$ represent positions of two masses, and the rest ones $\{x_t(3), x_t(4)\}$ stand for velocities. Parameters of this system are set as $K=1,\ m_1 = 0.5$, and $m_2 = 2$. The system is affected by unbounded disturbance $w_t \in \mathbb{R}$, whose data are generated based on an auto-regressive time series model.

The control objective is to steer the system towards the origin from the initial position $x_0 = [0.2\ 1\ -0.1\ 0.1]^{\rm T}$, and coefficient matrices in the quadratic objective function are specified as $Q=5I_4,\ R = 1$, and $Q_f = I_4$. At each stage $t$, we consider chance constraints imposed upon velocities of two masses:
\begin{equation}
\mathbb{P} \left \{ -0.38 \le x_t{(3)} \le 0.38,\ -0.38 \le x_t{(4)} \le 0.38 \right \} \ge 1 - \epsilon_t,\ t=1,\cdots,H,
\label{eq:statecc}
\end{equation}
along with hard constrains on control inputs:
\begin{equation}
u_t \in F_u = \left \{ u \left | -1.6 \le u \le 1.6 \right . \right \}, \ t=1,\cdots,H.
\label{eq:velocitycc}
\end{equation}
The saturated DF with the hyperbolic tangent function is adopted.

\subsection{Comparisons on Control Performance and Sampling Complexity}
First, we report control performances of the two-mass-spring system under different horizon lengths. For comparisons, the control decisions derived by the proposed DRMPC, classic SSMPC \cite{calafiore2006scenario} and RMPC \cite{zhang2013stochastic} are implemented in a RHC scheme, all of which provide appropriate probabilistic guarantee for chance constraints (\ref{eq:statecc}). Parameters in the probabilistic guarantee are set as $\epsilon_t \equiv 0.05$, $\beta_t \equiv 0.05$ for all stages, leading to different sample sizes required by various methods to safeguard the probabilistic guarantee. After formulating optimization problems as QPs, the solver {\tt cvx} in MATLAB \cite{grant2008cvx} is adopted and all computations are performed on a desktop computer with Intel Core i7-6700 processor at 3.40 GHz and 32 GB of RAM.

We perform simulations within 10 s ($1 \le t \le 100$), where the average cost-to-go $\frac{1}{T}\sum_{t=0}^{T-1}(x_t^{\rm T}Qx_t + u_t^{\rm T}Ru_t)$ as an empirical measure of control performance and the average CPU time are calculated and summarized in Table \ref{tab:ctrlperformance}. It can be observed that control performances obtained by different methods get improved with the length $H$ of control horizon increasing. The RMPC method \cite{zhang2013stochastic} appears to yield the worst performance, mainly because the hyper-rectangle uncertainty set leads to over-conservative solutions; however, it achieves the least computational time since the induced robust optimization problems are quite easy to solve. In other words, control performance is somewhat sacrificed for computational convenience by RMPC \cite{zhang2013stochastic}. The generic SSMPC method enjoys better control performance than RMPC \cite{zhang2013stochastic}, but the scale of optimization problem increases rapidly as the control horizon length $H$ and sample size increase. In practice, one could utilize parallel computing techniques to alleviate the computational burden. In contrast, the merits of the proposed DRMPC method can be demonstrated from various aspects. On one hand, the best control performance is achieved by DRMPC, mainly because fewer scenarios are theoretically required to reach the probabilistic guarantee, and the SVC-based uncertainty set can enclose scenarios non-conservatively. On the other hand, the computational costs are still acceptable thanks to an RO reformulation. Therefore, a better trade-off between control performance and computational burden is obtained by DRMPC.
\begin{table}[htbp]
\caption{Control Results of Two-Mass-Spring System $(\epsilon_t = 0.05, \beta_t = 0.05)$}
\centering
\begin{tabular}{c c c c c}
\hline
& & Avg. & Avg. & Sampling \\
& & Cost-To-Go & CPU Time (s) & Complexity \\
\hline
      & SSMPC & 5.7117 & 2.2020 & $N = 434$ \\
$H=5$ & RMPC & 6.0306 & 0.4664 & $N = 311$ \\
      & DRMPC & 5.2010 & 1.0029 & $N_{\rm train} / N_{\rm calib} = 300 / 59$ \\
\hline
      & SSMPC & 4.8119 & 5.0475 & $N = 577$ \\
$H=6$ & RMPC & 5.4727 & 0.4934 & $N = 361$ \\
      & DRMPC & 4.1595 & 1.4268 & $N_{\rm train} / N_{\rm calib} = 300 / 59$ \\
\hline
      & SSMPC & 4.3150 & 10.3079 & $N = 740$ \\
$H=7$ & RMPC & 5.3097 & 0.5709 & $N = 410$ \\
      & DRMPC & 3.9448 & 2.0197 & $N_{\rm train} / N_{\rm calib} = 300 / 59$ \\
\hline
\end{tabular}
\label{tab:ctrlperformance}
\end{table}

\begin{figure}[ht]
\centering
\includegraphics[width=0.96\textwidth]{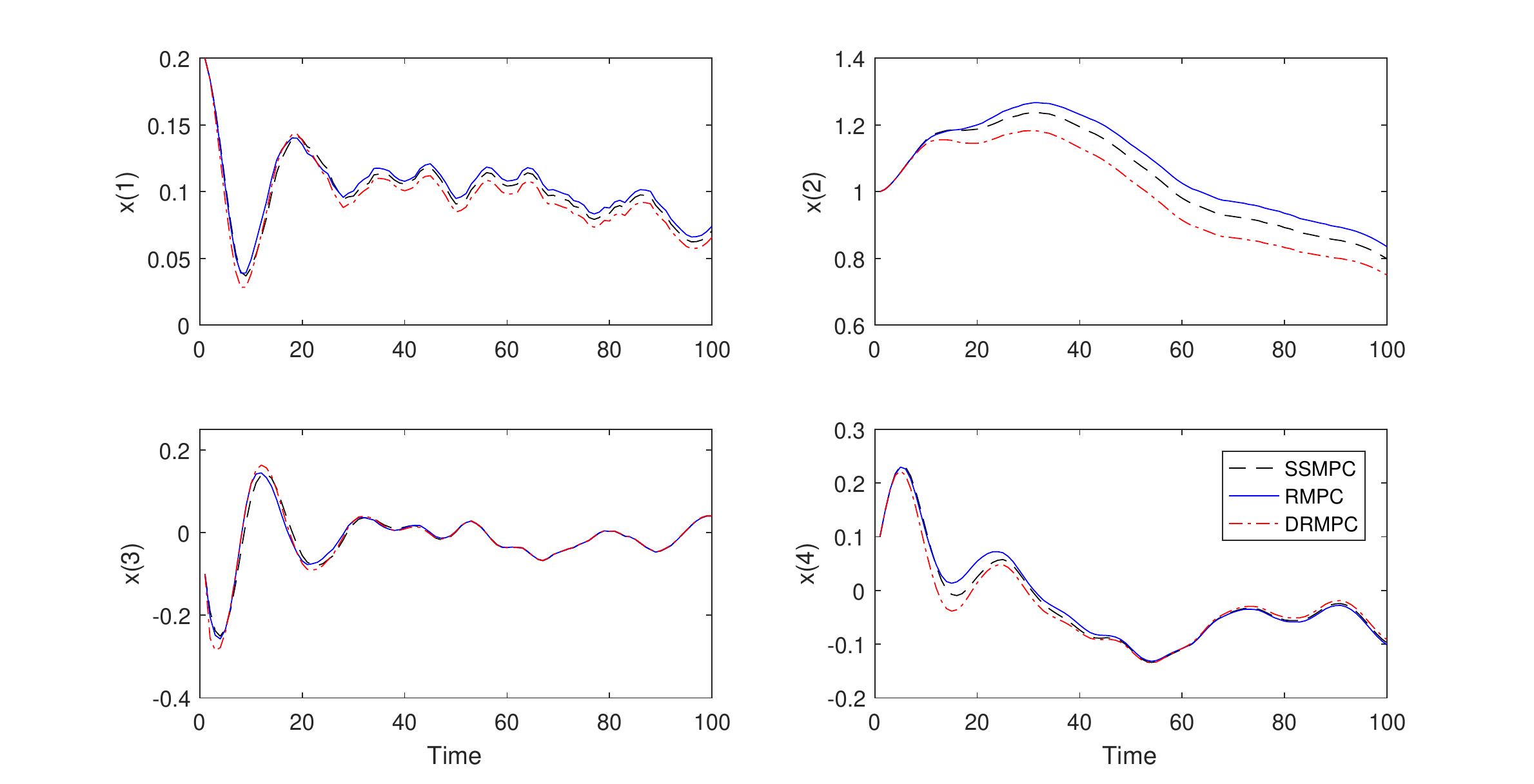}\\
\caption{Profiles of states in the two-mass-spring system with horizon length $H = 5$.}
\label{fig:stateprofile_H5}
\end{figure}

\begin{figure}[ht]
\centering
\includegraphics[width=0.5\textwidth]{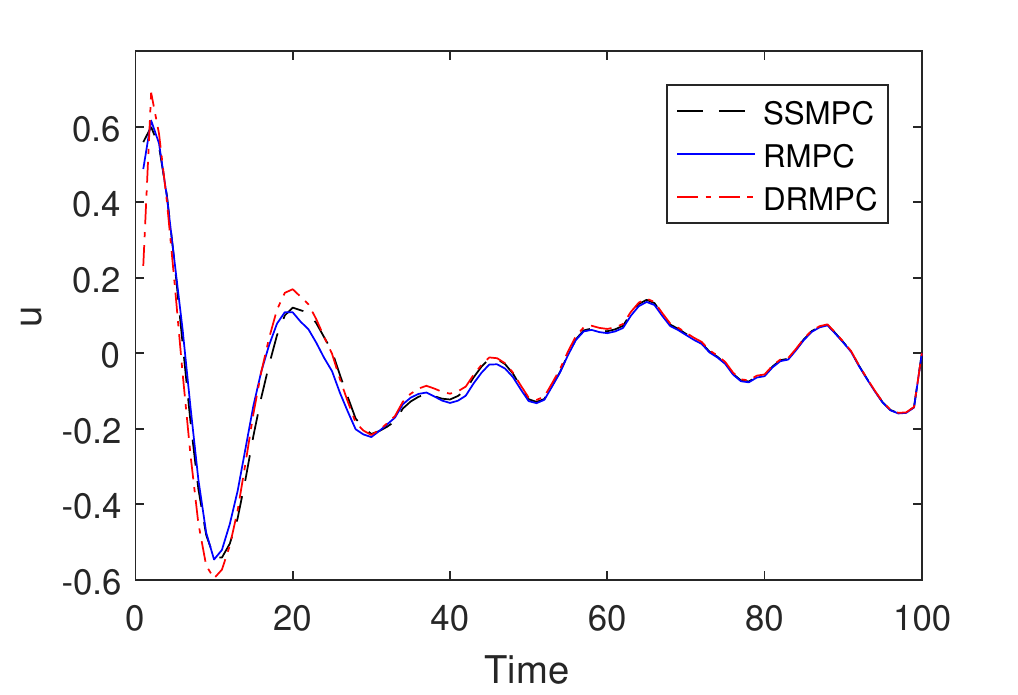}\\
\caption{Profiles of control input in the two-mass-spring system with horizon length $H = 5$.}
\label{fig:inputprofile_H5}
\end{figure}

\begin{figure}[ht]
\centering
\includegraphics[width=0.96\textwidth]{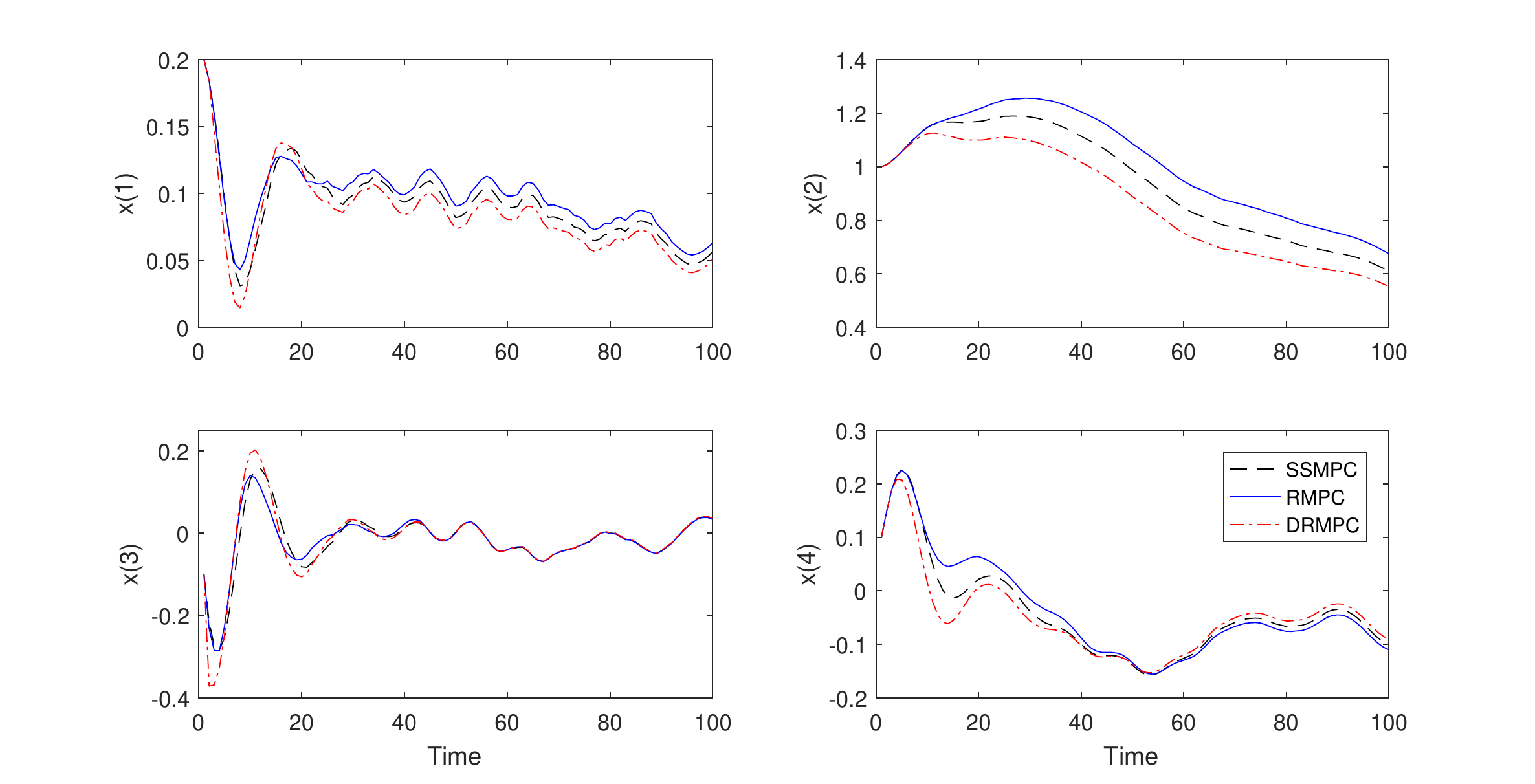}\\
\caption{Profiles of states in the two-mass-spring system with horizon length $H = 6$.}
\label{fig:stateprofile_H6}
\end{figure}

\begin{figure}[ht]
\centering
\includegraphics[width=0.5\textwidth]{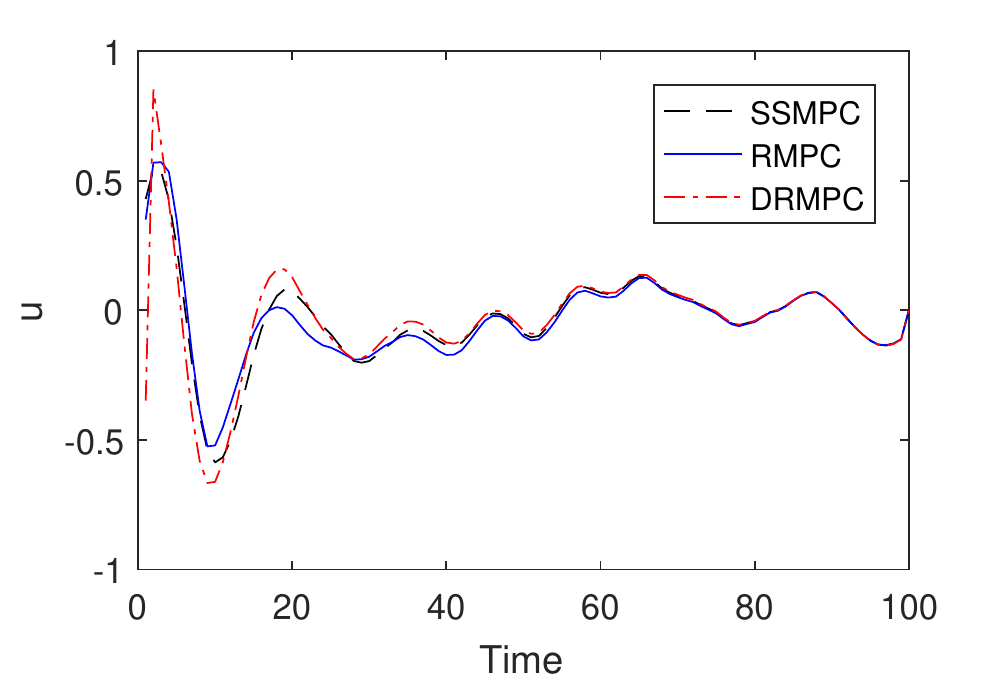}\\
\caption{Profiles of control input in the two-mass-spring system with horizon length $H = 6$.}
\label{fig:inputprofile_H6}
\end{figure}

\begin{figure}[ht]
\centering
\includegraphics[width=0.96\textwidth]{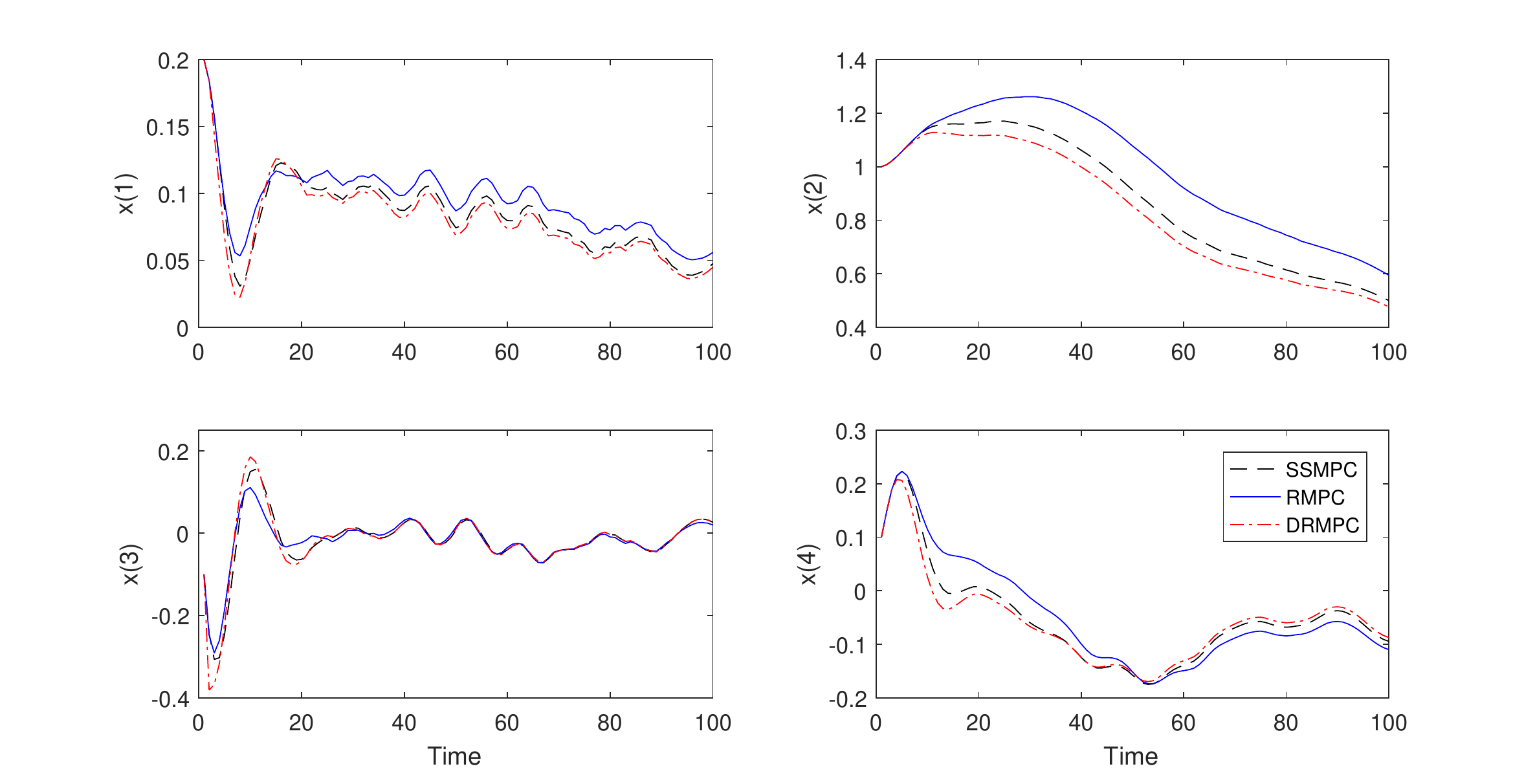}\\
\caption{Profiles of states in the two-mass-spring system with horizon length $H = 7$.}
\label{fig:stateprofile_H7}
\end{figure}

\begin{figure}[ht]
\centering
\includegraphics[width=0.5\textwidth]{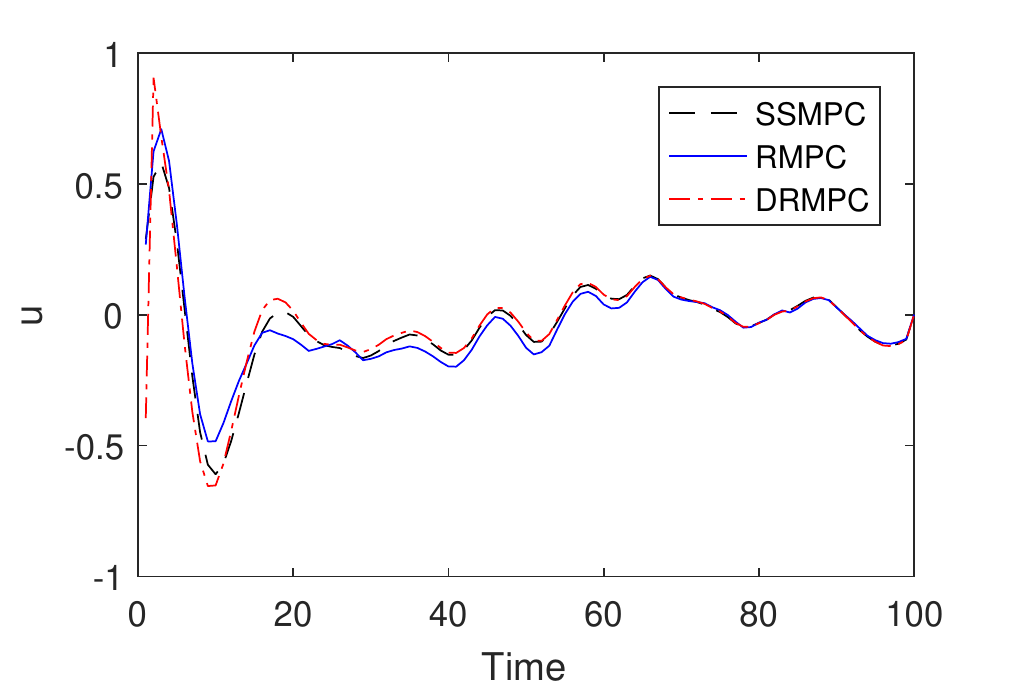}\\
\caption{Profiles of control input in the two-mass-spring system with horizon length $H = 7$.}
\label{fig:inputprofile_H7}
\end{figure}

To obtain more insights into the control performance, profiles of states and inputs under different horizon lengths $H$ are shown in Figures \ref{fig:stateprofile_H5} - \ref{fig:inputprofile_H7}. Due to the significant inertia and the influence of disturbance, a large settling time is typically required to steer states towards zero. In DRMPC $x(1)$ and $x(2)$ are closer to the origin during the transient process, especially the position of the second mass $x(2)$. Most importantly, the velocity $x(3)$ induced by DRMPC decreases more abruptly than SSMPC and RMPC, and gets closer to the pre-specified lower bound $-0.38$, which takes effect in driving masses towards the origin. This fundamentally owes to more aggressive control actions by DRMPC, as shown in Figures \ref{fig:inputprofile_H5}, \ref{fig:inputprofile_H6} and \ref{fig:inputprofile_H7}. It indicates that by effectively utilizing information within disturbance data, DRMPC yields compactly enclosing uncertainty sets and less conservative control actions, while still safeguarding the prior probabilistic guarantee with a smaller number of scenarios used.

\begin{figure}[ht]
  \centering
  \includegraphics[width=0.5\textwidth]{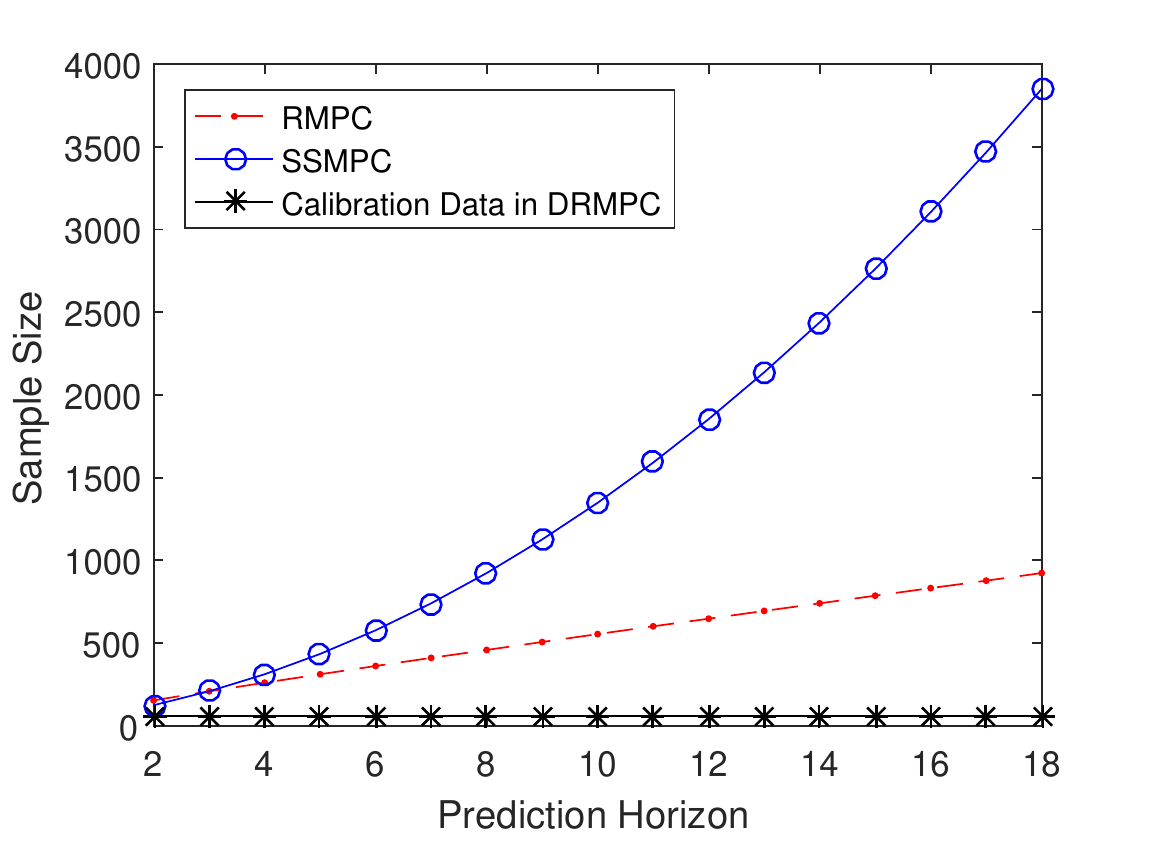}\\
  \caption{Sample sizes required for the probabilistic guarantee with different lengths of prediction horizon ($\epsilon = 0.05, \beta = 0.05$). }\label{fig:1}
\end{figure}

Next, we investigate the sample sizes required by different methods to attain the same probabilistic guarantee. In Table \ref{tab:ctrlperformance}, the sampling complexities required by SSMPC and RMPC increase with the horizon length $H$, while for DRMPC the theoretical number of calibration samples remains unchanged. We further vary the length of control horizon $H$, and more results are shown in Figure \ref{fig:1}, where the sample size computed by SSMPC explodes quadratically with $H$. This is because in SSMPC \cite{calafiore2006scenario} the number of decision variables induced by DF control policy can be calculated as $d_H = Hn_u + n_u n_w (H-1)H/2$. The method presented by \cite{zhang2013stochastic} achieves a sample size that grows approximately linearly with $H$, since the number of decision variables for learning a hyper-rectangle from scenarios is $2Hn_w$, which is proportional to $H$. By contrast, the proposed DRMPC method requires very few samples for calibration, which is a constant unrelated to $H$. Notice that a theoretical sample size for training data in DRMPC cannot be established, and the total sample size required by DRMPC is $N_{\rm train} + N_{\rm calib}$. In spite of this, because DRMPC typically entails a much smaller number of calibration samples $N_{\rm calib}$, there exists much room to choose adequate training samples $N_{\rm train}$ to derive a compactly enclosing uncertainty set, while still attaining a reduced total sampling complexity $N_{\rm train} + N_{\rm calib}$ in comparison with the other two methods.

\begin{figure}[ht]
  \centering
  \includegraphics[width=0.5\textwidth]{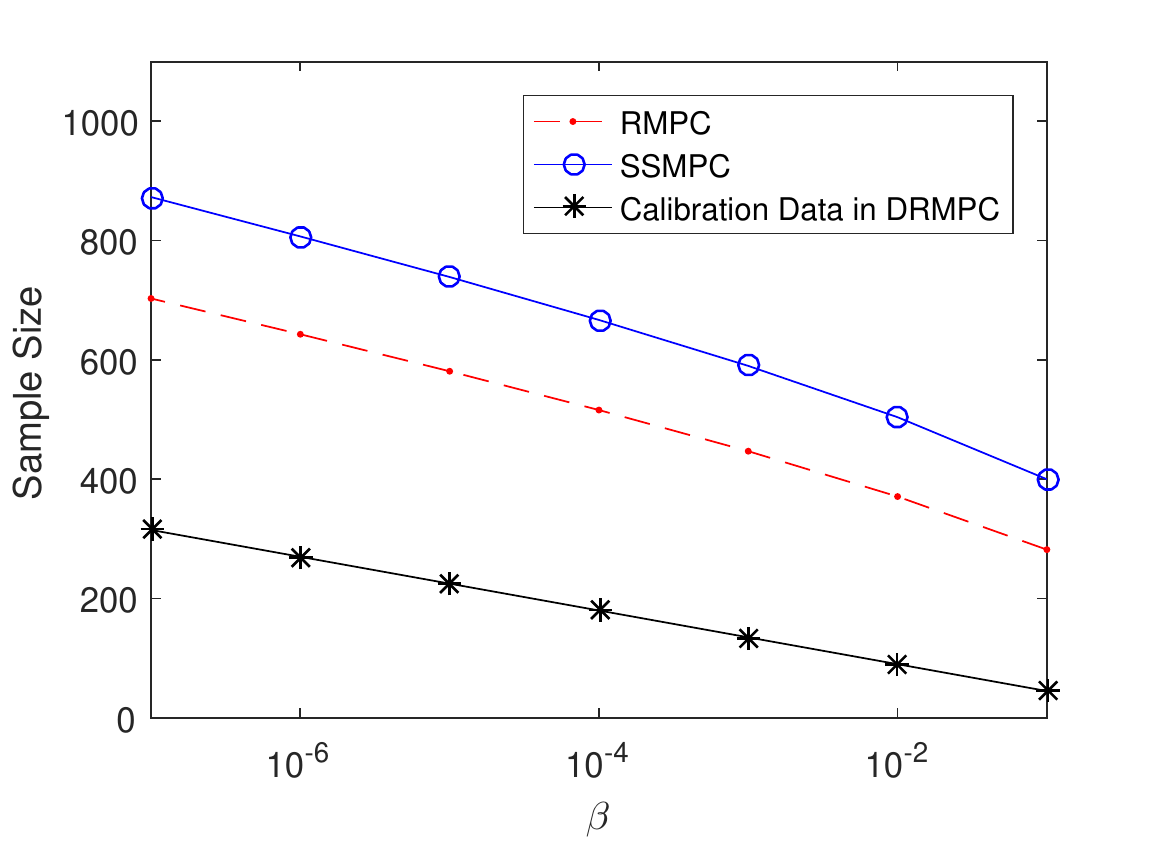}\\
  \caption{Sample sizes required for the probabilistic guarantee under different confidence levels ($\epsilon = 0.05, H = 5$). }\label{fig:2}
\end{figure}

\begin{figure}[ht]
  \centering
  \includegraphics[width=0.5\textwidth]{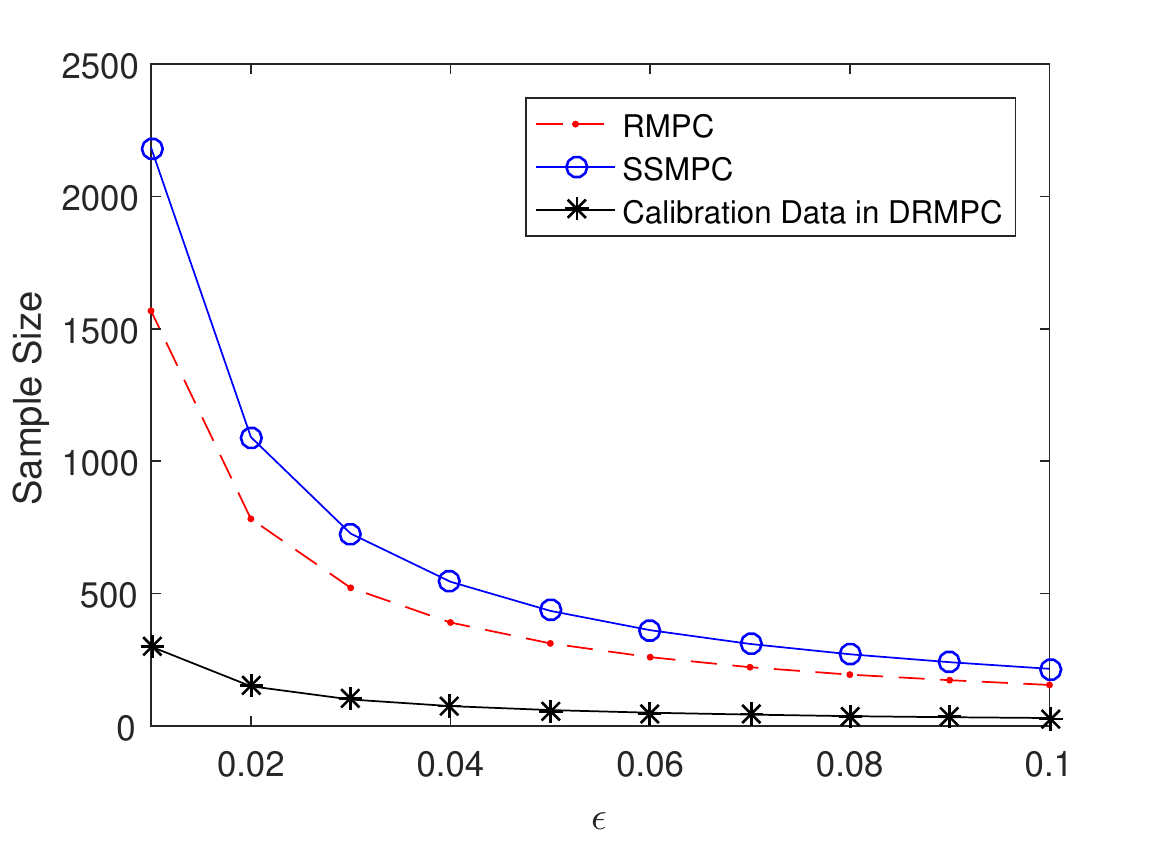}\\
  \caption{Sample sizes required for the probabilistic guarantee under different violation probabilities ($\beta = 0.05, H = 5$). }\label{fig:3}
\end{figure}

Next we investigate the effects of $\epsilon$ and $\beta$ on the sample size, as displayed in Figures \ref{fig:2} and \ref{fig:3}. We can observe that with values of $\epsilon$ and $\beta$ decreasing, the numbers of scenarios required by SSMPC \cite{calafiore2006scenario} and RMPC \cite{zhang2013stochastic} grow rapidly. It implies that traditional SSMPC \cite{calafiore2006scenario} and RMPC techniques in \cite{zhang2013stochastic} will gradually lose their practicability if a low violation probability $\epsilon$ or a high confidence level $1-\beta$ is pursued. By contrast, the sample size in the calibration phase required by DRMPC appears to be more practically acceptable with small $\epsilon$ or $\beta$. With a moderate sample size adopted in the training phase, pronounced advantages of DRMPC over classical approaches can be obtained.

\subsection{Comparisons on Empirical Probabilistic Guarantee}
The proposed DRMPC and RMPC \cite{zhang2013stochastic} share some similarities in that performance guarantee (\ref{eq:probgrt}) is first established on some uncertainty set $\mathcal{W}(\mathcal{D})$, namely the one based on SVC and the other one based on hyper-rectangle, which further ensures the performance guarantee (\ref{eq:probguarantee}) on final solutions to robust optimization problems. Before closing this section, we carry out Monte Carlo simulations to testify the probabilistic guarantees (\ref{eq:probgrt}) achieved by DRMPC and RMPC \cite{zhang2013stochastic}. Parameters are set as $\epsilon_t \equiv 0.05$, $\beta_t \equiv 0.05$ and $H = 10$. We execute 5,000 Monte Carlo runs, and in each run two different uncertainty sets are built based on random scenarios, and then the empirical violation probabilities $\epsilon$ of both the hyper-rectangle set and the SVC-based set are evaluated based on 20,000 independent data samples. By summarizing results in 5,000 runs, empirical confidence levels $\hat{\beta}$ in (\ref{eq:probgrt}) can be estimated, and the results are given in Table \ref{tab:montecarlo}. It can be observed that an over-conservative estimate of $\beta$ is obtained by RMPC \cite{zhang2013stochastic}, primarily because there is excessive coverage in the hyper-rectangle uncertainty set. By contrast, results of DRMPC are much less conservative, where different numbers of training samples are adopted $N_{\rm train}$ but the total sample size is always smaller than that in RMPC \cite{zhang2013stochastic}. In all four cases, all empirical estimates of $\beta$ are very close to the prespecified value 0.05. Another observation is that the number of training samples in DRMPC has minor effects on the probabilistic guarantee, since it only determines the initial ``shape" of the uncertainty set.
\begin{table}[htbp]
\caption{Statistics of Monte Carlo Simulations}
\centering
\begin{tabular}{c c c c c}
\hline
& Sample Size & $\hat{\beta}$ & mean of $\epsilon$ & var. of $\epsilon$ \\
\hline
RMPC & $N = 311$ & 0.0130 & 0.0257 & $8.16\times 10^{-5}$ \\
DRMPC & $N_{\rm train} / N_{\rm calib} = 100/ 59$ & 0.0486 & 0.0168 & 0.0027 \\
DRMPC & $N_{\rm train} / N_{\rm calib} = 150/ 59$ & 0.0492 & 0.0170 & 0.0028 \\
DRMPC & $N_{\rm train} / N_{\rm calib} = 200/ 59$ & 0.0492 & 0.0169 & 0.0028 \\
DRMPC & $N_{\rm train} / N_{\rm calib} = 250/ 59$ & 0.0496 & 0.0168 & 0.0027 \\
\hline
\end{tabular}
\label{tab:montecarlo}
\end{table}

\section{Case Study on Building Energy Control}
\begin{figure}[h]
\centering
\subfigure[Appearance]{
\includegraphics[width = 0.4\textwidth]{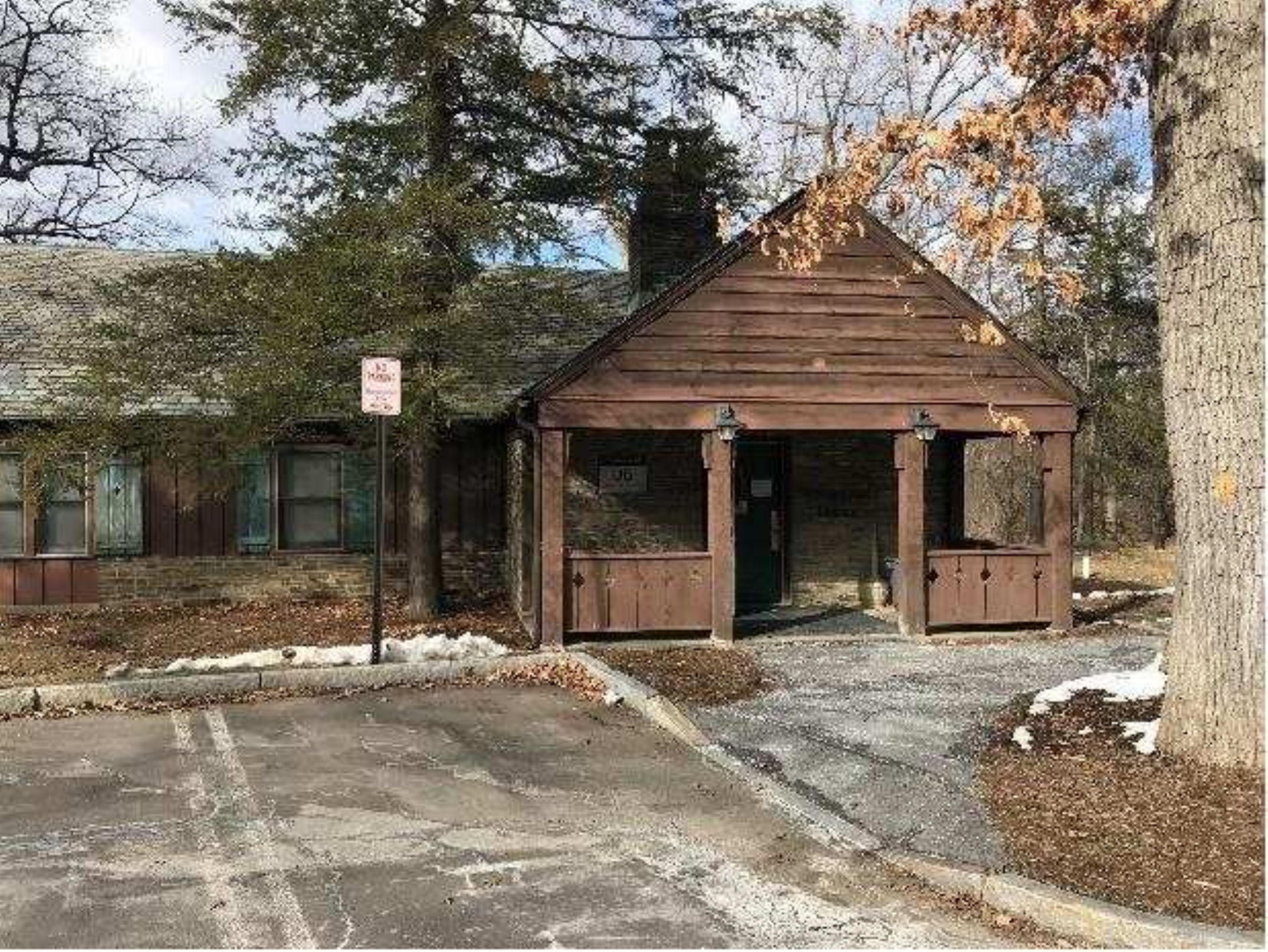}}
\subfigure[3-D modeling in BRCM toolbox]{
\includegraphics[width = 0.5\textwidth]{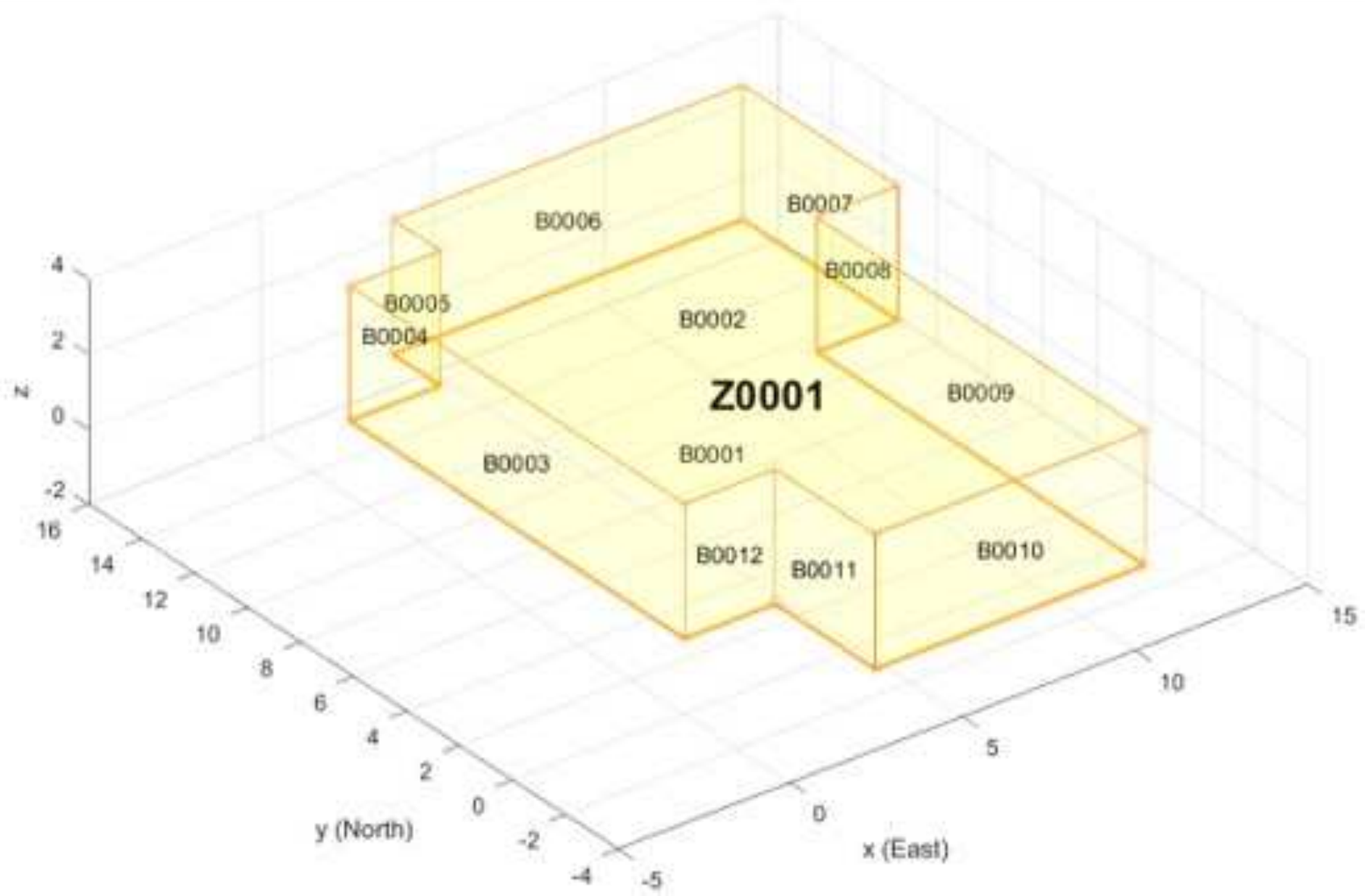}}
\caption{Toboggan Lodge at Cornell University, Ithaca, New York.}
\label{fig:lodge}
\end{figure}

Buildings are inevitably influenced by a variety of disturbances, such as outdoor temperature, solar radiation, etc., and the consideration of disturbances can enhance the efficiency of energy usage \cite{ma2009model,afram2014theory}. In this section, we aim to investigate the efficacy of the proposed DRMPC approach based on a simulated building control problem. It arises from the Toboggan Lodge, a one-floor campus building at Cornell University that has an irregular structure with multiple walls, as shown in Figure \ref{fig:lodge}. The Building Resistance-Capacitance Modeling (BRCM) MATLAB Toolbox developed by ETH Z{\"u}rich \cite{sturzenegger2014brcm} is adopted in this study for convenient generation of state-space models for buildings based on construction data. The building dynamics can be succinctly modelled as:
\begin{equation}
\begin{split}
x_{t+1} = Ax_t + B_u u_t + B_w w_t + B_v v_t
\end{split}
\end{equation}
where $x_t = [x_t(1)\ x_t(2)\ x_t(3)\ x_t(4)]^{\rm T}$ consists of room temperature, wall temperature, roof temperature, and floor temperature. There is only one input $u_t(1)$, namely the internal heating, and two deterministic external disturbances $v_t = [v_t(1)\ v_t(2)]^{\rm T}$, where $v_t(1)$ is the future prediction of external temperature, and $v_t(2)$ is the underground temperature, which is set as the annual average outdoor temperature \cite{zheng1993daily}. The model considered here is slightly different from (\ref{eq:statespace}), since forecasted values can be adopted to construct $v_t$, and the uncertainty in prediction error constitutes the additive disturbance $w_t$. We assume that the measurement of underground temperature is perfectly known, while the uncertainty in outdoor temperature forecasts is to be hedged against in building energy control. The step size of building simulation in BRCM is set as 1 h, thereby resulting in the following system matrices:
\begin{equation}
A = \left [
\begin{array}{c c c c}
0.0167 & 0.0048 & 0.1245 & 0.1409 \\
0.0005 & 0.0002 & 0.0039 & 0.0044 \\
0.0253 & 0.0073 & 0.3321 & 0.0617 \\
0.0244 & 0.0070 & 0.0526 & 0.3456 \\
\end{array}\right ],
\end{equation}
\begin{equation}
B_u = \left [
\begin{array}{c}
0.0986 \\
0.0029 \\
0.0288 \\
0.0275 \\
\end{array}
\right ],\ B_w = \left [
\begin{array}{c}
0.2536 \\
0.0070 \\
0.4450 \\
0.4477 \\
\end{array}
\right ],\ B_v = \left [
\begin{array}{c c}
0.2536 & 0.4596 \\
0.0070 & 0.9840 \\
0.4450 & 0.1287 \\
0.4477 & 0.1225 \\
\end{array}
\right ].
\end{equation}

The control objective of this building is to minimize the energy consumption while maintaining the room temperature within a comfortable range under weather prediction uncertainty. In the control objective, we set $Q = Q_f = 0$ and $R = I$, which penalizes energy consumption for heating. Time-varying thermal comfort constraints are imposed on the room temperature, which are formulated as chance constraints (\ref{eq:cc}). To be more specific, the room temperature $x_1(t)$ must be not lower than $21^{\rm o}{\rm C}$ during on-peak occupied hours (7am - 6pm), while in the rest off-peak hours, the room temperature $x_1(t)$ must be not lower than $15^{\rm o}{\rm C}$ to avoid unnecessary energy usage. This allows us to readily construct the matrix $F_t $ and vector $f_t$ in the chance constraint (\ref{eq:cc}). Limitations of heating power are modelled as hard constraints $0 \le u_t \le 80$ on control inputs. Then the original SMPC problem can be clearly constructed.

\subsection{Data Acquisition and Processing}
Historical predictions and real measurements of outdoor temperature at Ithaca, New York are available from Iowa Environmental Mesonet at Iowa State University.\footnote{https://mesonet.agron.iastate.edu/} Figure \ref{fig:temptrend_jan} exemplifies the profile of real measurements and 24-hour-ahead predictions in a certain day, from which prediction errors can be observed. In this way, a single scenario $\bold{w}$ representing prediction error uncertainty can be obtained by calculating the difference between measurements and predictions. Realistic weather data at Ithaca, New York, USA from Feb. 2016 to May 2017 have been collected. Data collected from Feb. 2016 to Dec. 2016 are used to construct scenarios of uncertainty, while data from Jan. 2017 to May 2017 are used for long-run simulations under various control strategies.
\begin{figure}[h]
\centering
\includegraphics[width = 0.55\textwidth]{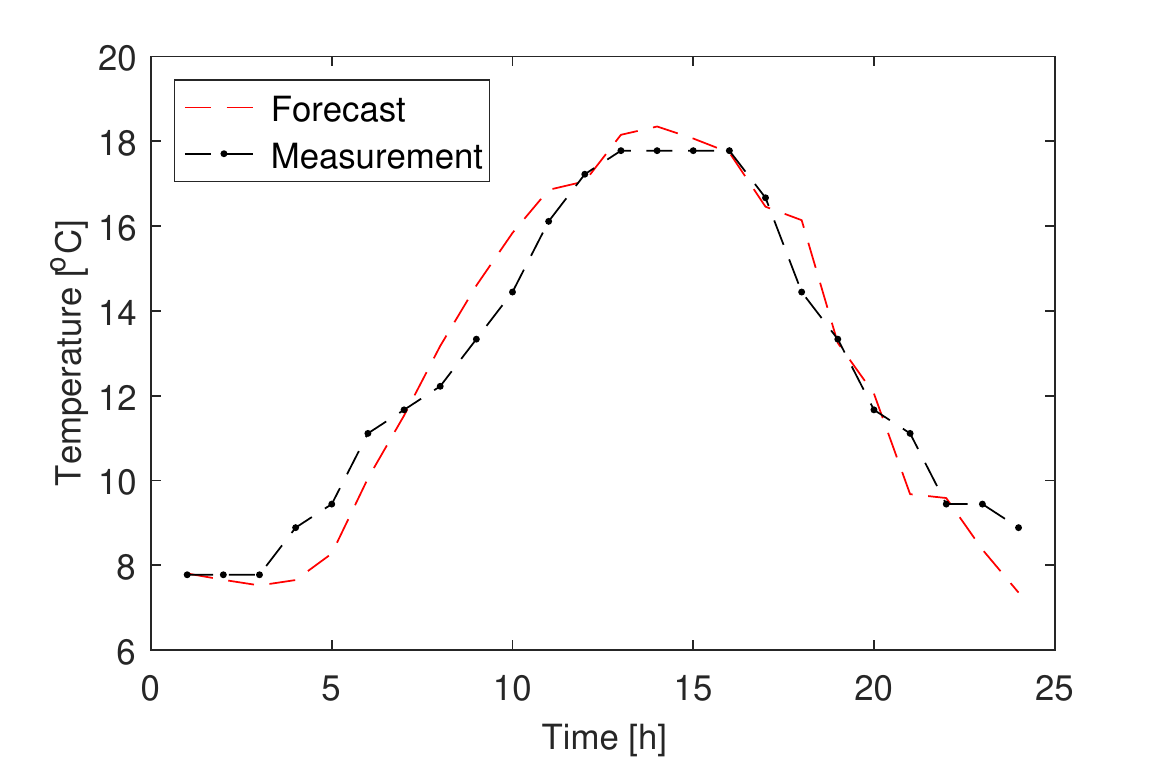}
\caption{24-hour-ahead forecast and real measurements of external temperature on June 9, 2017 at Ithaca, New York, USA.}
\label{fig:temptrend_jan}
\end{figure}

Here three control strategies under uncertainty are adopted, namely generic SSMPC \cite{calafiore2006scenario}, RMPC \cite{zhang2013stochastic} and the proposed DRMPC method, in a receding horizon fashion. The control horizon is set as $H = 5$, and the maximal violation probabilities and the confidence levels are set as $\epsilon_t = 0.05$ and $\beta_t = 0.10$ for all stages. Accordingly, 400 and 282 data samples are required by generic SSMPC \cite{calafiore2006scenario} and RMPC \cite{zhang2013stochastic}, respectively. As for the proposed method, $N_{\rm calib} = 45$ data samples are entailed in the calibration dataset $\mathcal{D}_{\rm calib}$, while $N_{\rm train} = 237$ scenarios constitute the training dataset $\mathcal{D}_{\rm train}$. In this way, under the same probabilistic guarantee, DRMPC has the same total sample size as RMPC \cite{zhang2013stochastic}, which is smaller than generic SSMPC \cite{calafiore2006scenario}.

\begin{figure}[h]
\centering
\includegraphics[width = 1.01\textwidth]{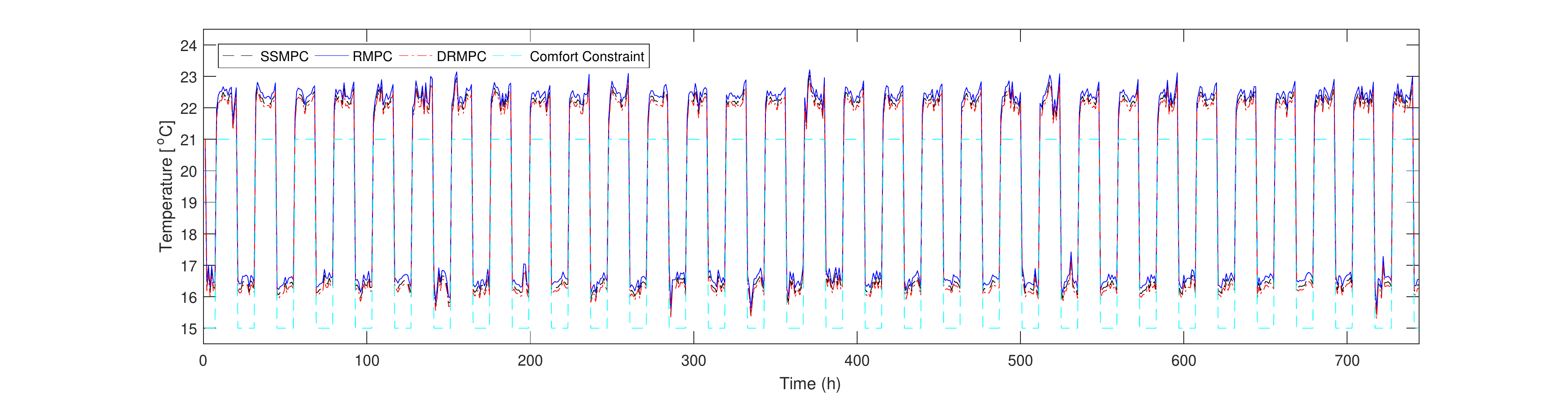}
\caption{Room temperature profile in Jan. 2017.}
\label{fig:temptrend_jan}
\end{figure}

\begin{figure}[h]
\centering
\includegraphics[width = 1.01\textwidth]{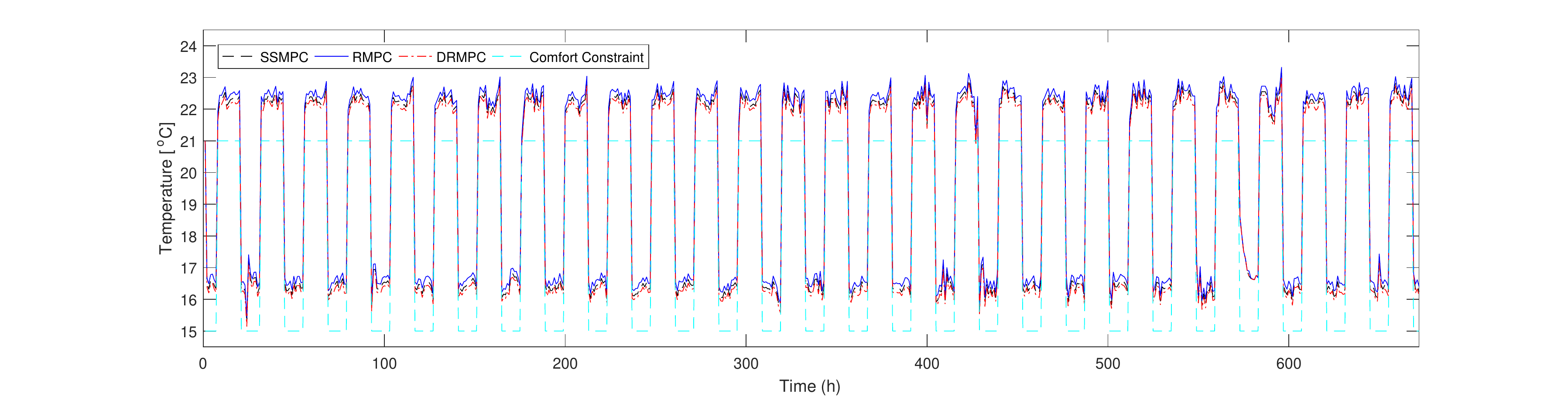}
\caption{Room temperature profile in Feb. 2017.}
\label{fig:temptrend_feb}
\end{figure}

\begin{figure}[h]
\centering
\includegraphics[width = 1.01\textwidth]{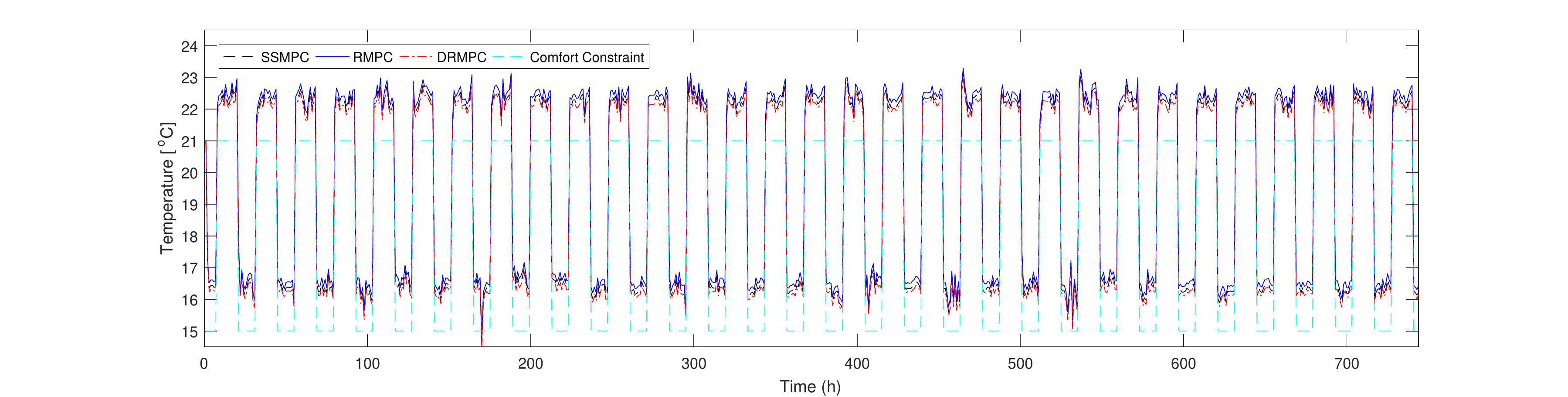}
\caption{Room temperature profile in Mar. 2017.}
\label{fig:temptrend_mar}
\end{figure}

\begin{figure}[h]
\centering
\includegraphics[width = 1.01\textwidth]{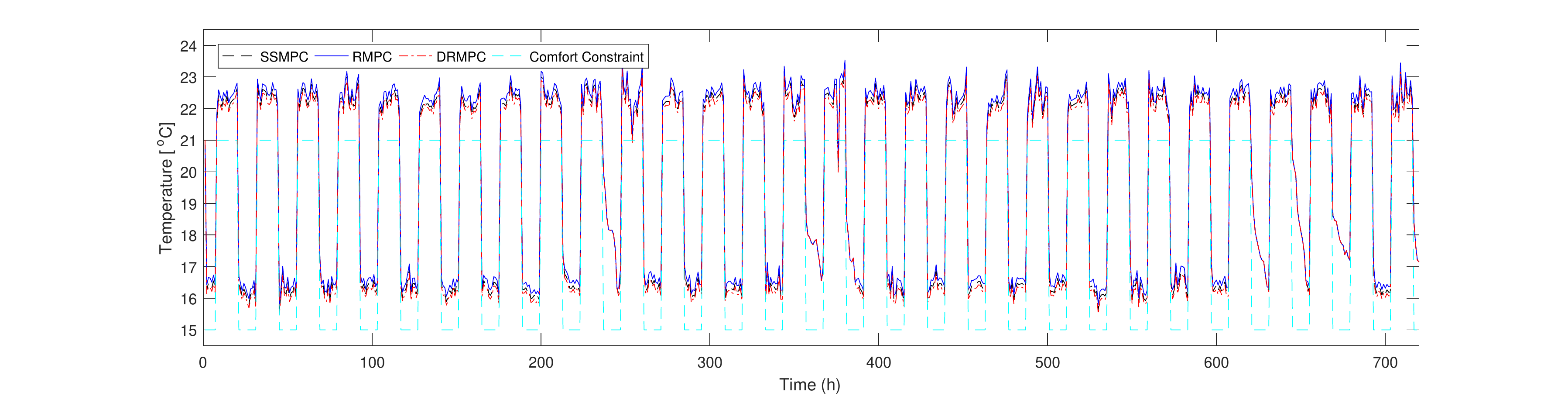}
\caption{Room temperature profile in Apr. 2017.}
\label{fig:temptrend_apr}
\end{figure}

\begin{figure}[h]
\centering
\includegraphics[width = 1.01\textwidth]{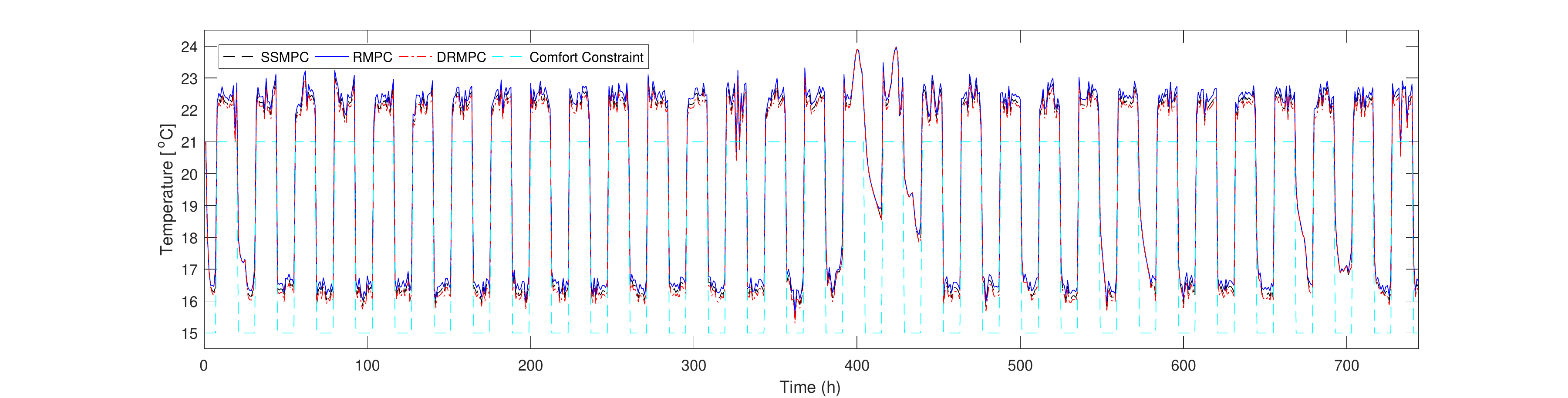}
\caption{Room temperature profile in May 2017.}
\label{fig:temptrend_may}
\end{figure}

\subsection{Control Results and Discussions}
Next, we investigate the long-term control performance of different strategies based on real measurements and predictions of outdoor temperature collected from Jan. 2017 to May. 2017, during which indoor heating is necessary. In Figures \ref{fig:temptrend_jan} - \ref{fig:temptrend_may}, profiles of room temperature in each month under different control strategies are depicted. It can be seen that, in all cases, the room temperature satisfies the thermal comfort constraint with some margins, which is a consequence of pursuing robustness. In particular, the proposed DRMPC provides the least conservative control performance, since room temperature is closest to the constraint boundary. This is also justified by Table \ref{tab:energyconsumption}, where the control performance in terms of monthly energy consumption is given. It can be seen that the proposed approach leads to the lowest energy usage, while RMPC \cite{zhang2013stochastic} leads to the highest. The reduction in conservatism by DRMPC is due to the fact that a small number of scenarios are required in the calibration phase, and the SVC-based uncertainty set achieves a compactly enclosing polytope and avoids unnecessary coverage of generic norm-based sets. Therefore, to attain the same probabilistic guarantee, less extreme realizations of uncertainty are to be protected against by DRMPC. In this way, the system can be operated closer to the constraint boundary in a less conservative manner, which eventually brings energy savings. In contrast, in RMPC one needs to ensure the robust constraints based on a hyper-rectangular uncertainty set, and in SSMPC a large number of scenarios are involved in the optimization problem, thereby inducing over-conservatism to various extent.

Another observation in Figures \ref{fig:temptrend_jan} - \ref{fig:temptrend_may} is that violations of comfort constraints may sometimes occur. Therefore, we calculate the empirical probabilities of comfort constraint violations in the long run from Jan. 2017 to May 2017, and the results are reported in Table \ref{tab:tempviolation}. Due to the inherent robustness of the RHC scheme, empirical probabilities in all cases are still far below the prespecified value $\epsilon_t \equiv 5\%$. However, the reduced conservatism of DRMPC can be seen from the results in April and May, where higher probabilities of comfort constraint violations are attained. This is because massive scenarios in SSMPC \cite{calafiore2006scenario} and the hyper-rectangular uncertainty set used in RMPC \cite{zhang2013stochastic} tend to involve more extreme realizations of weather prediction error, which in turn highlights the advantage of DRMPC in deriving less conservative solutions to chance-constrained optimal control problems.

\begin{table}[htbp]
\caption{Monthly Energy Consumptions for Temperature Control of Toboggan Lodge in 2017}
\centering
\begin{tabular}{c c c c c c}
\hline
& Jan. & Feb. & Mar. & Apr. & May \\
\hline
SSMPC [W/m$^2$] & 68,346 & 55,453 & 67,617 & 32,264 & 27,092 \\
RMPC [W/m$^2$] & 69,480 & 56,460 & 68,747 & 33,262 & 27,497 \\
DRMPC [W/m$^2$] & 67,591 & 54,785 & 66,865 & 31,604 & 25,916 \\
\hline
\end{tabular}
\label{tab:energyconsumption}
\end{table}

\begin{table}[htbp]
\caption{Empirical Probabilities of Comfort Constraints Violation of Building Climate Control in 2017}
\centering
\begin{tabular}{c c c c c c}
\hline
& Jan. & Feb. & Mar. & Apr. & May \\
\hline
SSMPC [\%] & 0 & 0.30 & 0.13 & 0.14 & 0.40 \\
RMPC [\%] & 0 & 0 & 0.13 & 0.14 & 0.27 \\
DRMPC [\%] & 0 & 0.30 & 0.13 & 0.28 & 0.54 \\
\hline
\end{tabular}
\label{tab:tempviolation}
\end{table}

\section{Conclusions}
In this paper we proposed a data-driven robust MPC approach, termed as DRMPC, to approximately solve scenario-based SMPC problems under uncertain disturbance. The idea is to first learn a high-density region as a data-driven uncertainty set from random realizations of disturbance and then solve a robust optimization problem to determine the control actions of the current time step and make predictions on future steps. The SVC-based uncertainty set was adopted to compactly capture the support of uncertainty in a data-driven manner. A novel calibration approach was proposed to further adjust the size of the uncertainty set and endow it with desirable probabilistic guarantees. We showed that the number of data samples required to attain the probabilistic guarantee can be further reduced in comparison with classical approaches. This is particularly beneficial for enhancing the applicability of SMPC. In addition, by constructing an compactly enclosing polytope of uncertainty, the conservatism of control actions can be greatly reduced. The issue of feasibility and stability of the proposed control scheme was also rigorously analyzed. Simulation results on a two-mass-spring system and a building energy control problem were presented, which well demonstrated the merits of the proposed DRMPC approach, that reduced conservatism and moderate computational burden can be obtained in order to reach the same level of probabilistic guarantee of SMPC.


\section*{Acknowledgments}
C. Shang acknowledges financial support from National Natural Science Foundation of China (Nos. 61673236, 61433001, 61873142). F. You acknowledges financial support from the National Science Foundation (NSF) CAREER Award (CBET-1643244). We also thank Haokai Huang, Wei-Han Chen, Jun Yang and Xiu Wang for their efforts in collecting weather data and developing the building control example.

\end{spacing}
\end{document}